\documentclass[msom,final]{informs-neutral}
\newif\ifarxiv
\arxivfalse

\newif\iflong
\longtrue

\newif\ifblind
\blindfalse

 \SingleSpacedXI

\usepackage{natbib}
\bibpunct[, ]{(}{)}{,}{a}{}{,}
\def\bibfont{\small}
\def\bibsep{\smallskipamount}

\usepackage[utf8]{inputenc}
\usepackage[T1]{fontenc}
\usepackage{microtype}

\usepackage{url}
\usepackage{hyperref}

\usepackage{booktabs}
\usepackage{longtable}

\usepackage{enumerate}
\usepackage{enumitem}
\usepackage{listings}

\usepackage{mathtools}
\usepackage{amsmath}
\usepackage{amssymb}
\usepackage{amsfonts}
\usepackage{nicefrac}
\usepackage{braket}
\usepackage{centernot}

\usepackage{float}
\usepackage{wrapfig}
\usepackage{graphicx}

\usepackage{ifthen}
\usepackage{color}
\usepackage{array}
\usepackage{doi}
\usepackage{setspace}

\newcommand{\N}{\mathcal{N}}
\newcommand{\T}{\mathcal{T}}
\newcommand{\B}{\mathcal{B}}
\newcommand{\G}{\mathcal{G}}
\newcommand{\LOS}{\mathcal{L}}

\TheoremsNumberedThrough
\ECRepeatTheorems
\EquationsNumberedThrough

\begin{document}

\MONTH{October}
\YEAR{2020}

\MANUSCRIPTNO{MSOM-20-502}

\TITLE{Optimal Resource and Demand Redistribution for Healthcare Systems Under Stress from COVID-19}
\RUNTITLE{Hospital Resource and Demand Redistribution for COVID-19}

\ifarxiv
\else
\ARTICLEAUTHORS{
    \AUTHOR{Felix Parker}
        \AFF{Department of Civil and Systems Engineering, The Center for Systems Science and Engineering, The Malone Center for Engineering in Healthcare, Johns Hopkins University, Baltimore, MD 21218, \EMAIL{fparker9@jhu.edu}}
    \AUTHOR{Hamilton Sawczuk}
        \AFF{Department of Applied Mathematics and Statistics, Johns Hopkins University, Baltimore, MD 21218, \EMAIL{hsawczu1@jhu.edu}}
    \AUTHOR{Fardin Ganjkhanloo, Farzin Ahmadi, Kimia Ghobadi}
        \AFF{Department of Civil and Systems Engineering, The Center for Systems Science and Engineering, The Malone Center for Engineering in Healthcare, Johns Hopkins University, Baltimore, MD 21218, \EMAIL{fganjkh1@jhu.edu}, \EMAIL{fahmadi1@jhu.edu}, \EMAIL{kimia@jhu.edu}}
}
\fi
\RUNAUTHOR{Parker et al.}

\ABSTRACT{
When facing an extreme stressor, such as the COVID-19 pandemic, healthcare systems typically respond reactively by creating surge capacity at facilities that are at or approaching their baseline capacity. However, creating individual capacity at each facility is not necessarily the optimal approach, and redistributing demand and critical resources between facilities can reduce the total required capacity. Data shows that this additional load was unevenly distributed between hospitals during the COVID-19 pandemic, requiring some to create surge capacity while nearby hospitals had unused capacity. Not only is this inefficient, but it also could lead to a decreased quality of care at over-capacity hospitals. In this work, we study the problem of finding optimal demand and resource transfers to minimize the required surge capacity and resource shortage during a period of heightened demand. We develop and analyze a series of linear and mixed-integer programming models that solve variants of the demand and resource redistribution problem. We additionally consider demand uncertainty and use robust optimization to ensure solution feasibility. We also incorporate a range of operational constraints and costs that decision-makers may need to consider when implementing such a scheme. Our models are validated retrospectively using COVID-19 hospitalization data from New Jersey, Texas, and Miami, yielding at least an 85\% reduction in required surge capacity relative to the observed outcome of each case. Results show that such solutions are operationally feasible and sufficiently robust against demand uncertainty. In summary, this work provides decision-makers in healthcare systems with a practical and flexible tool to reduce the surge capacity necessary to properly care for patients in cases when some facilities are over capacity.
}

\maketitle

\section{Introduction} \label{section:introduction}

Since the first confirmed case of a SARS-CoV-2 infection in Wuhan, China, COVID-19 has rapidly evolved into a global pandemic \citep{hui2020continuing}. As of September 2020, there have been more than 48 million reported cases worldwide with the United States reporting the most confirmed COVID-19 cases among all countries at nearly nine million. The world has seen more than one million and 236 thousand deaths \citep{dong2020interactive} as the result of the pandemic. COVID-19 has created an enormous burden on the capacities of healthcare systems across the globe, especially in Intensive Care Units (ICUs) \citep{covid2020forecasting}. These healthcare capacity concerns have prompted state and local governments to intervene by instituting widespread closures and stay-at-home orders \citep{mervosh_2020}. Yet as the number of cases continues to rise around the globe and in the United States (US), healthcare capacity remains a concern at the hot spots of the pandemic. 

In this paper, we present a framework to match the demand and available resources in healthcare networks. We, in particular, focus on scenarios in which individual hospitals in the network are under capacity stress, such as the COVID-19 pandemic. In this framework, we develop several allocation optimization models to achieve a system-level balance of loads by redistribution of patients (as opposed to resources such as ventilators or nurses) among different hospitals. The optimization models aim to minimize patient overflow in each hospital while also considering the operational constraints of patients transfers. We tested our proposed models on three case studies in New Jersey, Florida, and Texas using publicly available data gathered from state and healthcare agencies. We are currently in the process of adapting our models and tools for implementation in a large academic healthcare system in the US as it prepares for a potential new wave of SARS-CoV-2 infection in fall 2020. Our models will be used for optimal patient transfers within this healthcare system as well as informing strategic decisions on COVID-19 capacity management at each hospital.

\emph{Motivation and Impact}. With the arrival of the pandemic to the US and the strain on healthcare capacity, the healthcare system responded by both reducing demand through canceling elective surgeries and discharging or discharging existing patients earlier and also increasing capacity through calling additional personnel and creating new COVID-19-suitable beds as part of the guidelines from the Centers for Disease Control and Prevention (CDC) \citep{cdc_manage_2020}. While these measures have proven effective, delays in elective care may result in adverse effects for patients. Additionally Increasing COVID-19 capacity is costly, slow, and may not be feasible at all healthcare centers. At the same time, healthcare centers have made efforts to optimize the use of complementary resources, including extending the use of Personal Protective Equipment (PPE) \citep{centers2020strategies}. However, in practice, these methods are often insufficient or unsustainable. Another compounding organizational drawback to conventional healthcare center responses is that often they are organized at individual healthcare centers level rather than system levels. Here, we aim to take a step toward a more coherent response from the system of healthcare entities in times of chronic demand surge.

On a local level, patients tend to choose healthcare centers based on reputation or distance, leading to unbalanced patient loads across healthcare centers, which in turn decreases the overall quality of patient care \citep{varkevisser2012patients, drevs2013patients}. A similar phenomenon can be observed on a larger scale where different regions experience varying pandemic-related demands at different times. Hence, while some healthcare centers may be operating at capacity, others might have additional capacity to spare. 
Current forecasts show that COVID-19 burdens healthcare centers far beyond their current capacity, especially in ICUs, which calls for a model that optimally distributes patients across healthcare centers to minimize overall capacity and resource shortages.
Addressing the issue of unbalanced patient loads requires considering healthcare centers at the system level -- across hospitals, counties, and even states -- where there is potential for capacity pooling to lead to more efficient resource use \citep{vanberkel2012efficiency, sims2015evidence, chod2005resource}. Such approaches have been employed in ad-hoc manners during the pandemic. For instance, \citep{leon2020using} describe patient transfer among a system of hospitals in New York despite operational challenges and \cite{cdc2020key} includes two instances of a situation where the states of Michigan and Texas opted to utilize such methods in extreme case surge intervals and were able to partially alleviate the situation. Given the ongoing capacity concerns, healthcare systems are increasingly considering system-level interventions and patient transfers to maximize the utilization of available resources.

\emph{Models and Case Studies.} 
In this work, we propose several optimization models to balance system loads in the face of an extreme stressor. Our models produce optimal and operationally feasible redistributions of newly admitted patients across healthcare systems to minimize total patient overflow and consequently reduce the strain on healthcare workers and infrastructure. Because we consider only newly admitted patients, the process is made easier for hospitals to decide on new patients. The transfer of the patients happens at their point-of-entry to the healthcare center, which is often the Emergency Department (ED). Although our reallocations are assumed to take place immediately upon patient arrival at a hospital, depending on EMS operational considerations, these transfers could potentially take place before arrival, rerouting patients on the way to hospitals. 

The proposed models encompass a larger scope of operational situations. The models developed in this work are capable of also redistributing nurses and other critical resources across healthcare systems simultaneously with patients, hence, further reducing shortage that may result from closed beds. We additionally consider various types of patients who might require potentially different care paths, e.g., acute-level care, ICU-level care, or any combination of them. This distinction plays an important operational role in capacity management as the available capacity, the patient demographic, and the LOS for each level of care is vastly different. These differences may become even more pronounced in subsequent waves of the pandemic. Additionally, we consider uncertainties that exist in capacity management at hospital and system level, especially during a pandemic, from the number of patients and the level of their required care based on their possibly uncertain LOS. To address such issues, we consider a range of forecasts and distributions and propose a robust optimization model. As it will be outlined in detail in subsequent sections, it is worth noting that our models directly use the output of external patient count forecasts, which allows for further flexibility in using forecasting models. 

Implementation of such a model provides important managerial insights. In a general setting, the recommendations from our models give care centers experiencing excess demand a data-driven plan to reduce their loads and maintain adequate service quality. The flexibility afforded by these load reductions will also compound, allowing decision-makers to respond better to changing system demands. When combined with accurate forecasts, insights from our models can indicate vulnerable centers where an advanced expansion of capacity would be beneficial to the system as a whole. Additionally, the results on specific case studies considered in this work show that, for example, it was possible to reduce patient overflow by nearly 90\% in New Jersey, one of the worst impacted states during the first wave of the pandemic. It is worth to mention that, due to a challenge in our case studies with access to publicly and trustworthy available data, we created a publicly available repository of relevant healthcare and forecasting data that can be used by other researchers and practitioners
\ifblind
\else
\hspace{-0.2em}. The URL for this site is: \url{https://jhu-covid-optimization.github.io/covid-data/}.
\fi

A summary of the contributions of this work is as follows:
\begin{enumerate}
    \item We provide a demand allocation optimization model to distribute demand across health care entities levels in the presence of extreme demand. This model allows for different patient care paths, each with an associated length of stay distribution for each bed type.
    \item We show how to incorporate additional penalties and constraints to encourage the operational feasibility of solutions in the proposed model. We also propose a model that is robust against the most prominent sources of uncertainty in this setting.
    \item We compile healthcare data relevant to our modeling and make it publicly available along with our model implementations and maintain an interactive website that allows users to specify model parameters and generate custom reallocation solutions.
    \item We present several case studies and our models are under implementation at a large academic health system in the US. 
\end{enumerate}

The rest of this paper is organized as follows. Related literature and the recent developments during the COVID-19 pandemic are presented in Section \ref{section:related}. Section \ref{section:methods} provides details on the proposed allocation models to improve load-balance in healthcare networks. Section \ref{section:data} discusses the details of the data utilized in this study. The results of our analysis for the models and two case studies are illustrated in Section \ref{section:results}. We conclude by discussing future directions in Section \ref{section:conclusion}.

\section{Related Work} \label{section:related}

In this section, we highlight some of the existing literature on topics related to this work. We start by introducing some of prominent papers on general resource allocation. Related literature on healthcare operations, modelings for pandemics, and pandemic response are also provided. Finally, we discuss the specific domain of healthcare resource allocation and works related to the COVID pandemic. 

The general problem of resource allocation has been well studied in the literature, leading to the development of many different optimization methods in diverse domains. Earlier works in this field, such as the work done by \cite{hegazy1999optimization}, show a parallel development of optimization and heuristic algorithms to tackle the general problem of resource allocation. Other works in the literature concentrate on combining traditional resource allocation optimization models with different approaches \citep{cruz2004particle}. More recently, the general problem of optimizing resource allocation has gained significant attention in complex fields such as networks and systems \citep{ren2018latency, kuchuk2015two, tychogiorgos2014optimization}. Additionally, there is an extensive body of literature discussing the issue of robustness in general optimization problems and also in resource allocation settings \citep{mulvey1995robust,najafi2013transportation}. Among such works, in the general optimization context, \cite{bertsimas2004price} introduce the notation of the price of robustness and robust formulations in forms of linear optimization problems and \cite{gabriel2014recent} discuss recent work in the field of robust optimization in general.

Optimization in the settings of healthcare applications has also been the subject of many previous works. Most related to ours is the study done by \cite{mills2020surge} which discusses approaches for deploying and optimizing surge capacity in hospitals. Other works in the literature tackle similar problems in healthcare. \citep{luscombe2016dynamic} consider the problem of real-time dynamic emergency department (ED) scheduling with the goal of minimizing patient wait times. Similarly, \cite{otegbeye2015designing} develop a tool which utilizes retrospective ED data to generate optimal nurse shift schedules. \cite{litvak2008managing} present mathematical support for a system-wide pool of reserved emergency ICU beds to improve patient care and successful admission rates efficiently through network-wide cooperation. However, most of such works discuss normal and anticipated workflow or mass casualty incidents which have short term effects on the hospitals.  However, there have also been some studies that consider extended high demand periods in their approaches. For instance,  \cite{toner2006hospitals} emphasize that in the presence of an event such as a pandemic, cooperation between different healthcare centers becomes crucial to accommodate extreme system-wide stress. Nonetheless, much of the existing literature on optimal healthcare operations only considers the setting of a single healthcare center or does not consider high-demand settings where some healthcare centers rapidly hit their load capacity. For extensive surveys regarding operations research applications in healthcare systems, the reader is referred to the works of \cite{papageorgiou1978some} and  \cite{rais2011operations}.

Additional works explore healthcare system responses to pandemics specifically. For example, \cite{brandeau2004allocating} explore the problem of general allocation of resources in response to a pandemic, while \cite{toner2006hospitals} develop a set of best practices for influenza pandemic preparedness and response, and \cite{halpern2020urge} \cite{cooper2006delaying} explore the challenges and limitations of increasing critical care supply and shutting down air travel respectively in the presence of a pandemic. Others  consider more granular levels of analysis. For instance, as a general response to the COVID pandemic, \cite{centers2020strategies} propose operational strategies to optimize the use of face-masks and other personal protective equipment and \cite{mehrotra2020rapidly} build a framework for rapidly moving patient primary care to a virtual modality. Some more recent studies have explored ethical approaches to rationing critical resource allocations during the COVID pandemic \citep{emanuel2020fair, white2020framework}. \cite{judson2020rapid} consider building a patient self-scheduling and self-triage tool designed to decrease patient wait times during the COVID pandemic.

Another important and relevant area of research is methods for tracking and modeling the progression and effects of pandemics. This topic has been central in the literature \citep{pei2018forecasting}, and the literature has grown considerably since the emergence of the COVID pandemic \citep{petropoulos2020forecasting, jewell2020caution, roda2020why, perc2020forecasting}. Recently, \cite{lewnard2020incidence} and \cite{rees2020covid} retrospectively analyzed and synthesized COVID healthcare data to provide key qualitative and quantitative insights into the pandemic's nature and spread. \cite{dong2020interactive} provide valuable data by tracking cases across all nations and in the United States at the county level, and \cite{weissman2020locally} and \cite{covid2020forecasting} provide forecasts for patient loads at the hospital and state level respectively. In general, the existence of such data and tools is paramount to developing meaningful pandemic operations optimization models.

We close the related works section by pointing out existing works on the problem of resource allocation using central stockpiles during high demand periods like pandemics. This problem has been previously studied by \cite{arora2010resource} and \cite{mehrotra2020model} in deterministic and stochastic settings respectively. More specifically, the issue of optimal vaccine allocation has been addressed by \cite{longini1978optimization} and \cite{lampariello2020effectively} in non-COVID and COVID settings respectively and operational considerations of COVID patient transfers are discussed in \cite{cdc2020key} and \cite{kain2020safe}. As a real world example of patient allocation schemes during the COVID pandemic, \cite{leon2020using} present the successes and challenges associated with patient transfers among New York City hospitals. In the most relevant existing literature, \cite{bai2014incentive, sun2014multi} and \cite{lacasa2020flexible} discuss the problem of optimal patient allocation in a pandemic. \cite{lacasa2020flexible} consider the problem of distributing a single resource or demand across a regular geometric graph with healthcare centers as vertices. They provide solutions as a set of resource transfers using random search optimization. However, it should be noted that such solutions are usually not guaranteed to be optimal, only correspond to a single time step, do not support complementary resources or secondary operational constraints, and are not robust against any particular data uncertainties. Other related works in the literature include the study done by \cite{bai2014incentive} that casts the problem of patient load-balancing in a pandemic setting as a max-flow problem, offering a solution containing a set of patient incentives to optimally balance healthcare center loads. However, rather than considering minimizing the capacity overflow, they seek to minimize patient cost, suggesting that their model may be more applicable in a lower-demand setting. Finally, \cite{sun2014multi} produce a model similar to the models in this work, however, similar to \cite{bai2014incentive}, they attempt to minimize patient travel distance and assume a constant length of stay for each group of patients. They also provide fewer secondary operational constraints to improve practical solution feasibility. Here, we extend the capabilities of the models in the literature by proposing a new modeling methodology for the optimal allocation of patients and/or complementary resources subject to a variety of secondary operational constraints and potential uncertainties in different parameters like the number of admitted patients, beds and supplementary resources.

\section{Methodology} \label{section:methods}

In this section we present a series of linear optimization (LP and MILP) models to solve the multi-period demand and resource redistribution problems.

The demand redistribution problem is as follows: given a set of nodes and time periods, with nominal demand at each node during each period and fixed capacity at each node, determine the optimal quantity of demand to transfer between each pair of nodes during each time period. The specific metric according to which a solution is optimal may depend on context in real world settings; however, in this work we consider two primary objectives: minimize the total amount of surge capacity which must be created or balance the normalized load between nodes each time period. The resource redistribution may be defined analogously, with resource transfers as the decision variable rather than demand transfers.

We study these problems specifically in the context of hospital systems facing an extreme stress, where nodes represent hospitals, time periods are days, and demand is the number of patients that need care. Capacity is determined by the number of hospital beds as well as critical resources such as nurses, PPE, and ventilators. However, the methods developed are not specific to this context and are applicable to other domains.

In Section \ref{section:methods:patients} we develop the demand redistribution models. In Section \ref{section:methods:combined} we develop the resource redistribution model, which can be used independently or in tandem with the demand redistribution. Finally, in Section \ref{section:methods:evaluation} we describe the metrics used to evaluate the solutions of each model.
The details of our assumptions, parameters, variables and other notation for all models are summarized in Table \ref{tab:modelinfo} at the end of this section.

\subsection{Demand Redistribution Models} \label{section:methods:patients}

In this section we develop the demand redistribution models, starting with the base model upon which the subsequent models are built. The first extension of the model incorporates multiple demand and capacity types with different properties. In the healthcare setting these correspond to patient groups with different care paths, and to different bed types respectively. The next extension adds a series of operational costs and constraints that are optional, but potentially valuable to decision-makers. Finally, the model is made robust against reasonable levels of demand uncertainty.

In order to define the scope of our solutions to the demand redistribution problem we have made some assumptions.
First, we assume that demand redistribution decisions are made by a centralized decision-maker and that each node will follow decisions exactly.
Second, we assume that the decision-maker has full access to required data on capacity and demand (potentially with some uncertainty).
Third, we assume that demand is transferred in the same time period as it arrives at a node. In the context of patient transfers this means we assume transfers are made after a patient arrives at a hospital, but before they are admitted. This could include patients who have just arrived, or patients who have had some stay in the ED.
We also assume that demand arrives and then decays with some function over time. In the case of patient redistribution this means that patients have a length of stay and that the number of patients remaining $t$ days after they were admitted is a function of $t$. We use one minus the cumulative distribution function of the estimated LOS distribution as the decay function for patients. The accuracy of this approximation increases with the number of admitted patients.
Finally, in the context of patient transfers, we assume that decisions on which patients to transfer can be made effectively by decision-makers at each node as these decisions are out of the scope of this work. This holds for any discrete, non-homogeneous demand.
Any secondary assumptions that have been made will be noted as they arise.

\subsubsection{Base Demand Redistribution Formulation.} \label{section:methods:patients:base}
We begin by formulating a simple linear optimization model to solve the demand redistribution problem.
\begin{subequations}
\begin{align}
	\underset{\omega}{\text{minimize}} & \quad \sum_{i \in \N} \sum_{t \in \T} \omega_{i,t} \label{eq:PA_obj}\\
	\text{subject to} 
	& \quad \sum_{j \in \N} s_{i,j,t} \leq p_{i,t}  & \forall i \in \N, t \in \T \label{base:con:sent}\\
	& \quad \alpha_{i,t} - b_i \leq \omega_{i,t}  & \forall i \in \N, t \in \T \label{base:con:overflow}\\
	& \quad \alpha_{i,t} = (p_{i,0} - \sum_{t'=1}^{t} d_{i,t'}) + \sum_{j \in \N} s_{i,j,t} \nonumber\\& \quad + \sum_{t'=1}^{t} \{ [1 - \LOS(t - t')] [p_{i,t'} + \sum_{j=1}^{N} (s_{j,i,t'} - s_{i,j,t'})]\} &\forall i \in \N, t \in \T \label{base:expr:alpha}\\
	& \quad s_{i,j,t} = 0 & \forall (i,j) \in \overline{E(G)}, t \in \T \label{base:con:graph}\\
	& \quad s_{i,j,t} \geq 0, \quad \omega_{i,t} \geq 0 & \forall i,j \in \N, t \in \T \label{base:con:s0}
\end{align}
\end{subequations}

In this formulation, the objective is to minimize the total surge capacity that must be created to accommodate the demand. Required surge capacity is referred to as overflow, and demand is referred to as patients.
Only newly admitted patients may be transferred between nodes, which is enforced by constraint \eqref{base:con:sent}, ensuring that the number of patients transferred away does not exceed the number of newly admitted patients $p_{i,t}$. The overflow for node $i$ at time $t$, denoted by $\omega_{i,t}$, is defined by constraint \eqref{base:con:overflow} in conjunction with the constraint that $\omega_{i,t} \geq 0$. It is also assumed that transfers between certain nodes can be infeasible (e.g., based on distance or existing relationships between hospitals). These relationships and restrictions can be formally modeled as a graph $G$ where all healthcare centers are considered as vertices, and there exists and edge from node $i$ to node $j$ if patients can be transferred from node $i$ to node $j$. This assumption is captured by constraint \eqref{base:con:graph}.

Expression \eqref{base:expr:alpha} represents the number of active patients at node $i$ at time $t$. The first term of this expression captures the number of remaining initial patients. The third term incorporates the cumulative patient length of stay distribution $\LOS$ to capture the number of admitted and reallocated patients remaining at node $i$ at time $t$.
Since patient transfer between care centers typically requires resources from both centers on the day of the transfer, we count the transferred patients as active at both nodes $i$ and $j$ at time $t$, which is captured by the second term of the expression.
More details about the parameters and variables of this formulation are provided in Table \ref{tab:modelinfo}.

Note that in this model we have elected to relax the problem so that the number of transfers may be a real number rather than constraining it to be integral. Integrality in the number of transfers is left as an optional constraint in Section \ref{section:methods:patients:optional}. The purpose of this relaxation is to enable the problem to be solved using linear programming rather than integer programming, allowing non-trivially sized problems can be solved efficiently. We justify this relaxation by assuming that the number of transfers is large so that the error caused by rounding is relatively small.

Overflow minimization is the primary objective considered in this work, and will be used in all subsequent demand redistribution models. However, as previously mentioned, it is not the only possible objective. Another objective that can be used is load balancing. Load is defined as demand normalized by capacity, so it measures how much stress a node is under and can be compared between nodes. Balancing load can be beneficial for a system because it means that the stress is distributed evenly. This objective will also use the capacity of the system efficiently, so there will not be excess capacity at one node while another node is over capacity. It is therefore implicitly performing overflow minimization as well.
However, load balancing typically involves more patient transfers than overflow minimization would.
The formulation for load balancing is as follows:
\begin{subequations}
\begin{align}
	\underset{s,\lambda}{\text{minimize}} & \quad \sum_{i \in \N} \sum_{t \in \T} \lambda_{i,t} \label{eq:load:obj} \\
	\text{subject to} 
	& \quad l_{i,t} = \frac{\alpha_{i,t}}{b_{i}} & \forall i \in \N, t \in \T \\
	& \quad l_{i,t} - \frac{1}{|\N|} \sum_{j \in \N} l_{j,t} \leq \lambda_{i,t} & \forall i \in \N, t \in \T \\
	& \quad -(l_{i,t} - \frac{1}{|\N|} \sum_{j \in \N} l_{j,t}) \leq \lambda_{i,t} & \forall i \in \N, t \in \T
\end{align}
\end{subequations}

These models assume a single care path for patients, and a single limiting capacity, total staffed beds. In reality however, patients have different care needs and will be assigned to different types of beds accordingly. Additionally, these needs will change over time, so patients may move between bed types.
In the next section, we expand the base model to take into account patients with different care paths and their varying resource requirements.

\subsubsection{Demand Redistribution Formulation with Demand and Capacity Types.} \label{section:methods:patients:group}

In this section, we extend the base LP model discussed in Section \ref{section:methods:patients:base} to provide a model capable of considering a setting where with a set of patient groups $\G$ and a set of bed types $\B$. Each patient is assigned to a specific patient group upon admission. During the course of treatment, a patient will follow a particular care path consisting of a sequence patient group, potentially with different associated bed types. Once again we refer to demand as patients here and resources as beds.

To model the care paths, directed graph $G_{group}$ is defined such that each node corresponds to a specific phase of treatment with its own resource type requirements, referred to as a patient group, and each edge connects a patient group to the following patient group along a particular care path. In this setting, $G_{group}$ is a assumed to be a disjoint union of directed trees with edges directed toward their own roots (also known as an in-forest graph). This assumption disallows cycles in a patient group transfer scheme and simultaneously requires each patient to have a unique transfer path until they are discharged. 

For formulation purposes, a function $f$ mapping $\G$ to $\B$ is defined, implying that each patient group is associated with exactly one bed type. 
\begin{subequations}
\begin{align}
	\underset{\omega}{\text{minimize}} & \quad \sum_{\beta \in \B} \sum_{i \in \N} \sum_{t \in \T} \omega_{\beta,i,t} \label{group:obj}\\
	\text{subject to}
	& \quad \sum_{j \in \N} s_{g,i,j,t} \leq p_{g,i,t}  &\forall g \in \G, i \in \N, t \in \T \label{group:con:sent} \\
	& \quad \sum_{g \in img\set{f^{-1}(\beta)}}(\alpha_{g,i,t} + \sum_{j \in \N} s_{g,i,j,t}) - b_{\beta,i} \leq \omega_{\beta,i,t}  &\forall \beta \in \B, i \in \N, t \in \T \label{group:con:overflow} \\
    & \quad \alpha_{g,i,t} = p_{g,i,0} + \sum_{t'=1}^t (\chi_{g,i,t'} - \gamma_{g,i,t'}) & \forall g \in \G, i \in \N, t \in \T \label{group:expr:alpha} \\
    & \quad \chi_{g,i,t} = p_{g,i,t} + \sum_{g' : g' \sim g} \gamma_{g',i,t} + \sum_{j \in \N}(s_{g,j,i,t}-s_{g,i,j,t}) & \forall g \in \G, i \in \N, t \in \T \label{group:expr:beta}\\
    & \quad \gamma_{g,i,t} = d_{g,i,t} + \sum_{t'=1}^t \ell_{g}(t-t')\chi_{g,i,t'} & \forall g \in \G, i \in \N, t \in \T \label{group:expr:gamma} \\
	& \quad s_{g,i,j,t} = 0  & \forall g \in \G, (i,j) \in \overline{E(G)}, t \in \T \label{group:con:graph} \\
	& \quad s_{g,i,j,t} \geq 0, \omega_{\beta,i,t} \geq 0 & \forall g \in \G, \beta \in \B, i,j \in \N, t \in \T \label{group:con:s0} 
\end{align}
\end{subequations}

In the above model, constraints \eqref{group:con:sent} through \eqref{group:con:s0} are generalizations of the constraints \eqref{base:con:sent} through \eqref{base:con:s0} from the base model introduced in Section \eqref{section:methods:patients:base}. 
Expression \eqref{group:expr:alpha} represents the number of active patients in group $g$ at node $i$ at time $t$. The first term captures initial patients, the second term captures the assumption that patients will take up resources in both the sending and receiving hospitals on the day they are transferred, and the third term accounts for the sum of net active patient changes in group $g$ at node $i$.
Expression \eqref{group:expr:beta} represents the number of patients entering group $g$ in node $i$ at time $t$. The first term captures admitted patients into group $g$, the second term sums patients leaving other groups that transfer to the group $g$, and the third term includes net patient transfers in group $g$ at node $i$ at time $t$.
Finally, expression \eqref{group:expr:gamma} represents the number of patients leaving group $g$ in node $i$ at time $t$.
The first term of Expression \eqref{group:expr:gamma} captures initial patients discharged while the second term calculates the number of patients leaving group $g$ for another group or being discharged, at node $i$ at time $t$.
Notice that although expressions \eqref{group:expr:alpha} through \eqref{group:expr:gamma} do not immediately constitute a closed form expression for the number of active patients, they provide a method for computing one recursively which is guaranteed to terminate by our requirement that $G_{\text{group}}$ be a in-forest. More details on the parameters and variables of this formulation are provided in Table \ref{tab:modelinfo}.

The care path group patient allocation model developed here is capable of minimizing patient overflow across multiple patient groups and bed types in a given system. In the following section we provide a collection of additional constraints and penalties that can be added to the so far developed scheme, which add flexibility and capability to the models to account for probable operational considerations in the optimal patient allocation problem.

\subsubsection{Optional Constraints and Considerations.} \label{section:methods:patients:optional}

We now build on the group patient LP formulation to add additional penalties and constraints to ensure the operational feasibility of solutions. In particular, we add penalties and constraints to avoid solutions to the model presented in Section \ref{section:methods:patients:group} which may be difficult or impossible for healthcare systems to follow. While these additions to the model are not specific to COVID patient transfers, they are inspired by this context, and will be described as such. Details on the additional variables and parameters are provided in Table \ref{tab:modelinfo}.

\begin{subequations}
In the following formulations, penalty represents a sum of terms that can be added to the objective function of the model.
\begin{align}
	& s_{g,i,j,t} \in \mathbb{N} & \forall g \in \G, i,j \in \N, t \in \T
\end{align}
The first additional constraint we consider is to ensure that every transfer is of integral size. Often, demand is an integer-valued quantity, such as patients, which cannot be divided, so it is important that this is an option in the model. However, it is optional because it converts the model into an integer program, which is far harder to solve computationally, and because the continuous relaxation may be a good approximation in many cases.
\begin{align}
	&\sum_{g \in img\set{f^{-1}(\beta)}} \alpha_{g,i,t} \leq \max{b_{\beta,i}} &\forall \beta \in \B, i \in \N, t \in \T, \alpha'_{\beta,i,t} \leq b_{\beta,i} \label{opt:con:shortage1}\\
	&\sum_{g \in img\set{f^{-1}(\beta)}} \alpha_{g,i,t} \leq \alpha'_{\beta,i,t} &\forall \beta \in \B, i \in \N, t \in \T, \alpha'_{\beta,i,t} > b_{\beta,i} \label{opt:con:shortage2}
\end{align}
Constraint $\eqref{opt:con:shortage1}$ ensures that no node experiences a patient overflow at a time when they otherwise would not, and constraint \eqref{opt:con:shortage2} ensures that no node overflow is made more severe by our reallocation. These constraints may be critical to get hospitals to buy-in to the solution in practice because they ensure that no hospital will have to make their situation worse beyond their capacity.
\begin{align}
	\text{penalty} &= C_{\text{sent}} \sum_{g \in \G} \sum_{i,j \in \N} \sum_{t \in \T} s_{g,i,j,t} \label{opt:penalty:sent}
\end{align}
This penalizes the model for the total number of patients it transfers. Keeping everything else constant, a solution which transfers fewer patients is desirable because it is easier to implement operationally. In practice, there are many optimal or near-optimal solutions to the base demand redistribution model, so adding this penalty with a small $C_{\text{sent}}$ helps the model select an optimal solution that is operationally sound.
\begin{align}
	& \sum_{g \in \G} \sum_{i,j \in \N} \sum_{t \in \T} s_{g,i,j,t} \leq S_{\text{total}} & \\
	& s_{g,i,j,t} \leq S_{\text{max}} & \forall g \in \G, i,j \in \N, t \in \T
\end{align}
We now constrain the total amount of patients transferred, as well as the size of an individual transfer. These have the same motivations as penalty \eqref{opt:penalty:sent}, but could be more useful if a system can determine upper bounds on how many patients they can actually transfer.
\begin{align}
	\text{penalty} &= C_{\text{smooth}}\sum_{i,j \in \N} \sum_{t \in \T} \delta_{i,j,t} \\
	&\delta_{i,j,t} \leq \sum_{g \in \G}(s_{g,i,j,t - 1} - s_{g,i,j,t}) && \forall i,j \in \N, t \in \T \setminus \set{1}
	\label{opt:con:smooth1} \\
	&-\delta_{i,j,t} \leq \sum_{g \in \G}(s_{g,i,j,t - 1} - s_{g,i,j,t}) && \forall i,j \in \N, t \in \T \setminus \set{1}
	\label{opt:con:smooth2}
\end{align}
This is a penalty on the absolute value of the difference between patient transfer quantities from node $i$ to node $j$ at time $t-1$ and time $t$. The corresponding variable $\delta_{i,j,t}$ is determined by constraints \eqref{opt:con:smooth1} and \eqref{opt:con:smooth2}. These constraints are motivated by the idea that a more consistent transfer rate from node $i$ to node $j$ may be more operationally feasible.
\begin{align}
	\text{penalty} &= C_{\text{setup}} \sum_{i,j \in \N} \rho_{i,j} \\
	&M \sum_{g \in \G} \sum_{t \in \T}(s_{g,i,j,t} + s_{g,j,i,t}) \geq \rho_{i,j} &\forall i,j \in \N, j > i \label{opt:con:setup1} \\
	&m \sum_{g \in \G} \sum_{t \in \T}(s_{g,i,j,t} + s_{g,j,i,t}) \leq \rho_{i,j} &\forall i,j \in \N, j > i \label{opt:con:setup2}\\
	&\rho_{i,j} \in \set{0,1} & \forall i,j \in \N, t \in \T \label{opt:con:rho01}
\end{align}
This creates a setup cost -- that is a penalty for each pair of hospitals that transfer between themselves, which represents the potential overhead cost off establishing a transfer relationship. Transferring patients to a fewer number of other hospitals would likely make transfers easier to manage for the origin hospital. Constraints \eqref{opt:con:setup1} and \eqref{opt:con:setup2} define the binary variable $\rho_{i,j}$ accordingly.
\begin{align}
	&M \sum_{g \in \G} \sum_{j \in \N} s_{g,i,j,t} \geq \nu_{k,i,t} &\forall k \in \set{1,2}, i \in \N, t \in \T \label{opt:con:sendreceive1} \\
	&m \sum_{g \in \G} \sum_{j \in \N} s_{g,i,j,t} \leq \nu_{k,i,t} &\forall k \in \set{1,2}, i \in \N, t \in \T \label{opt:con:sendreceive2} \\
	&\nu_{k,i,t} + m \sum_{g \in \G} \sum_{j \in \N} \sum_{t' = t}^{\text{min}\set{t+T_{\text{switch}}, T}} s_{g,j,i,t'} \leq 1 &\forall k \in \set{1,2}, i \in \N, t \in \T \label{opt:con:sendreceive3} \\
	&\nu_{1,i,t}, \nu_{2,i,t} \in \set{0,1} & \forall i,j \in \N, t \in \T  \label{opt:con:nu01}
\end{align}
Constraints \eqref{opt:con:sendreceive1} through \eqref{opt:con:sendreceive3} enforce a minimum gap of $T_{\text{switch}}$ days between when a node may send and receive patients. This is done by introducing the binary variables $\nu_{1,i,t}$ and $\nu_{2,i,t}$ and using them to constrain the sent and received patients for each window of $T_{\text{switch}}$ days. In practice it does not make sense to rapidly switch between sending and receiving patients, or to send and receive at the same time, so these constraints ensure that does not happen.
\begin{align}
	\text{penalty} &= C_{\text{distance}} \sum_{g \in \G} \sum_{i,j \in \N} \sum_{t \in \T} d_{i,j} s_{g,i,j,t}
\end{align}
Here we penalize the total distance traveled per patient (or unit of demand). There can be a monetary cost, as well as an increase in risk to the patient's health, that is proportional to the time it takes to transfer them. It therefore can make sense to try to balance the time that patient transfers take with the primary objective.
\begin{align}
	&\sum_{g \in \G} s_{g,i,j,t} \in \set{0} \cup [S_{\text{min}}, S_{\text{max}}), \quad \delta_{i,j,t}, \phi_{\beta,i,t} \geq 0 & \forall \beta \in \B, i,j \in \N, t \in \T \label{opt:con:short1} 
\end{align}
Next, constraint $\eqref{opt:con:short1}$ makes $s$ a semi-continuous variable, meaning that it can either be equal to zero, or be in the interval $[S_{\text{min}}, S_{\text{max}})$. The effect of this is to set upper and lower bounds on the size of a transfer. Note that setting a lower bound on the size of the transfer makes the problem non-convex and therefore involves an integer variable.
\begin{align}
    &\alpha'_{\beta,i,t} = \sum_{g \in img\set{f^{-1}(\beta)}} [ p_{g,i,0} + \sum_{t'=1}^t (\chi'_{g,i,t'} -  \gamma'_{g,i,t'})] & \forall \beta \in \B, i \in \N, t \in \T \label{opt:expr:alpha'} \\
    &\chi'_{g,i,t} = p_{g,i,t} + \sum_{g' : g' \sim g} \gamma'_{g',i,t} & \forall g \in \G, i \in \N, t \in \T \label{opt:expr:beta'}\\
    &\gamma'_{g,i,t} = d_{g,i,t} + \sum_{t'=1}^t \ell_{g}(t-t')\chi'_{g,i,t'} & \forall g \in \G, i \in \N, t \in \T \label{opt:expr:gamma'}
\end{align}
Expressions \eqref{opt:expr:alpha'} through \eqref{opt:expr:gamma'} are analogous to expressions \eqref{group:expr:alpha} through \eqref{group:expr:gamma}, representing active patients, entering patients, and exiting patients for each group, except that these new expressions do not consider patient redistribution and $\alpha'_{\beta,i,t}$ is aggregated by bed type. Therefore expression $\eqref{opt:expr:alpha'}$ gives the number of active patients in bed type $\beta$ at node $i$ at time $t$ without any patient reallocation. This expression is used in some of the above constraints.
\end{subequations}

So far we have developed a demand redistribution scheme, with considerations for healthcare systems including operational constraints and different care paths and bed types. In the next section this scheme is extended to perform resource redistribution along with demand redistribution.

\subsection{Combined Demand and Resource Redistribution Model} \label{section:methods:combined}

In this section we build on the group patient LP formulation from section \ref{section:methods:patients:group} to allocate nurses as a primary transferable resource, along with the patients. In order for a patient to receive proper treatment we require that they both have a bed and adequate nurse care. First, a nurse allocation LP model is introduced. Secondly, a set of nurse specific optional constraints are added, addressing the issue of artificial shortage of static supply for nurses (compared to patients' dynamic demand). It is noteworthy that to simultaneously allocate patients and nurses, we can combine the variables and constraints from the models in this section with the variables and constraints from Section \ref{section:methods:patients} and take our new objective function to be a weighted sum of all of the individual model objective functions.

\subsubsection{Resource Redistribution Formulation.}
In the previous formulations we assumed that the resource availability for COVID patients is such that the capacity at each hospital is constant over time. This makes sense for some resources, such as hospital beds, which are reusable but cannot be moved between locations.
However, these are not good assumptions for other types of resources such as PPE or nurses. Both are able to be transferred and PPE cannot be reused indefinitely.
We therefore introduce a model which determines optimal resource transfers. This resource redistribution model can be used independently of demand redistribution, but its primary purpose is to be combined with the demand redistribution models.
\begin{subequations}
\begin{align}
    \underset{\theta}{\text{minimize}} & \quad \sum_{i \in \N} \sum_{t \in \T} \theta_{i,t} \label{resource:obj}\\
	\text{subject to} & \quad \sum_{j \in \N} \sigma_{i,j,t} \leq \eta_{i,t} & \forall i \in \N, t \in \T \label{resource:con:sent} \\
	& \quad q_{i,t} - \eta_{i,t} \leq \theta_{i,t} & \forall i \in \N, t \in \T \label{resource:con:shortage} \\
    & \quad \eta_{i,t} = n_{i} + \sum_{j=1}^{N} \sum_{t'=ts(t)}^{t} (\sigma_{j,i,t'} - \sigma_{i,j,t'}) &\forall i \in \N, t \in \T \label{resource:expr:eta} \\
    & \quad q_{i,t} = \sum_{\beta \in \B} \sum_{g \in img\set{f^{-1}(\beta)}} Q_{\beta} \alpha_{g,i,t} &\forall i \in \N, t \in \T \label{resource:expr:q} \\
	& \quad \sigma_{i,j,t} = 0 & \forall (i,j) \in \overline{E(G)}, t \in \T \label{resource:con:graph} \\
	& \quad \sigma_{i,j,t} \geq 0, \theta_{i,j} \geq 0 & \forall i,j \in \N, t \in \T \label{resource:con:s0}
\end{align}
\end{subequations}

Equations \eqref{resource:con:sent} through \eqref{resource:con:s0} are analogous to equations \eqref{base:con:sent} through \eqref{base:con:s0} in section \eqref{section:methods:patients}. Expression \eqref{resource:expr:eta} represents the resource supply and expression \eqref{resource:expr:q} represents the resource demand, which is taken to be proportional to the number of active patients. When distributing patients and resources simultaneously, $\alpha_{g,i,t}$ will depend on the model's patient allocation. 
Note that constraint \eqref{resource:expr:eta} assumes a reusable resource, of which there are a fixed quantity in the system. If the resource in question is not reusable and there can be supply from outside of the system we can replace constraint \eqref{resource:expr:eta} with the following:
\begin{align}
	&\eta_{i,t} = n_{i} + \sum_{t'=1}^{t} n_{i,t'} + \sum_{j=1}^{N} \sum_{t'=1}^{t} (\sigma_{j,i,t'} - \sigma_{i,j,t'}) && \forall i \in \N, t \in \T
	\label{resource:expr:eta_alt}
\end{align}

\subsubsection{Optional Costs and Constraints for Resource Redistribution.}
In this section we provide a set of optional constraints that can be used with the resource redistribution model in the case of a reusable resource. Specifically, the constraints provided account for the fact that if a node is experiencing a reusable resource shortage then that node should not have transferred resources to other nodes.
\begin{subequations}
\begin{align}
	&n_{i} \leq \eta_{i,t} &\forall i \in \N, t \in T, q_{i,t} \geq n_{i,t} \label{nurse:con:nurse1} \\
	&m (q_{i,t} - \eta_{i,t}) \leq \kappa_{i,t} &\forall i \in \N, t \in \T \label{nurse:con:kappa1} \\
	&1 + m (q_{i,t} - \eta_{i,t}) \geq \kappa_{i,t} &\forall i \in \N, t \in \T \label{nurse:con:kappa2} \\
	&\eta_{i,t} \geq n_i &\forall i \in N, t \in T, \kappa_{i,t} = 1 \label{nurse:con:nurse2} \\
    &\kappa_{i,t} \in \set{0,1} &\forall i \in \N, t \in \T \label{nurse:con:kappa0}
\end{align}
\end{subequations}

Constraints \eqref{nurse:con:kappa1}, \eqref{nurse:con:nurse2} and \eqref{nurse:con:kappa0} enforce that for any day that a node has a resource shortage, it has at least its initial supply. Optional penalties and constraints from \ref{section:methods:patients:group} can also be applied to the model, but we forgo rewriting them here for simplicity.

\subsection{Robust Optimization Model} \label{section:methods:robust}

Thus far our models have all implicitly relied on the assumption that we know the input data with certainty. In practice, however, regardless of whether our methodology is implemented prospectively or retrospectively, there is almost surely some extent of uncertainty in the data. We therefore construct a Robust Optimization model which ensures the optimal solution remains feasible for all scenarios in the uncertainty set we consider.
The model we present here addresses uncertainty on the number of admitted patients. The number of admitted patients is particularly impactful on the solution because it governs the number of active patients that need to be accommodated and it limits the number of patients that can be transferred. We also expect that this input will be the most uncertain data for prospective studies since we must rely on forecasting to generate the predictions.

We start with a box uncertainty set, meaning that we consider an upper and lower bound on the number of admitted patients, and assume that the true number of admitted patients must fall somewhere within this range. Specifically, $p_{i,t} \in P_{i,t} = [\bar{p}_{i,t} - \tilde{p}^{-}_{i,t}, \bar{p}_{i,t} + \tilde{p}^{+}_{i,t}]$ $\forall i \in \N, t \in \T$. We define $\bar{p}_{i,t}$ to be the nominal value of the number of admitted patients, which is what is used in the preceding models. Additionally, we assume that $\tilde{p}^{-}_{i,t}, \tilde{p}^{+}_{i,t} \geq 0$ and $\tilde{p}^{-}_{i,t} \leq \bar{p}_{i,t}$. Note that a distribution over the number of admitted patients is not specified. Rather, this method robustifies against any distribution with a fixed upper and lower bound, making it highly flexible.

Previous work in Robust Optimization, however, including \cite{bertsimas2004price} have shown that such models can be overly conservative with a high "price of robustness" -- that is they suffer a large increase to the optimal objective function value as compared with the certain optimization model. We therefore adopt an "uncertainty budget" inspired by \cite{bertsimas2004price}. Instead of assuming that the uncertain values can be anywhere between the upper and lower bound we constrain the uncertain values to be equal to the nominal value except on a parameter $\Gamma$ number of days, when it deviates to the upper or lower bound. Formally, $p_{i,t} = \bar{p}_{i,t} + \min \set{0,\xi_{i,t} \tilde{p}^{-}_{i,t}} + \max \set{0,\xi_{i,t} \tilde{p}^{+}_{i,t}}$ such that $\xi_i \in \{0,1\}^t, ||\vec{\xi}_i||_{1} = \Gamma$ $ \forall i \in \N$. This formulation captures the intuition that the true value of the number of admitted patients is generally equal to the nominal value, but will deviate on some days. The practitioner may vary $0 \leq \Gamma \leq |\N|$ where at the extremes the certain and fully robust are recovered respectively. Note that in practice we can relax the constraint $\xi_i \in \{0,1\}^t$ to $||\vec{\xi}_i||_{\infty} \leq 1$ and our solution method will recover the original constraint. Note also that due to the uncertainty budget, the numbers of admitted patients each day for a given node are dependent on one another. Consequently, a realization from the uncertainty set for a given node must consist of a sequence of admitted patients per day for that node.

A motivating factor for employing robust optimization and this uncertainty set in particular is that the robust model may be formulated as an LP with the same number of constraints and variables as the equivalent non-robust model, meaning that it can be solved efficiently in theory and practice. The details of this re-formulation may be found in the online supplement to this paper.

The following is the base patient allocation model presented in Section \ref{section:methods:patients:base} made robust against our identified uncertainty set.
\begin{subequations}
\begin{align}
	\underset{\omega}{\text{minimize}} & \quad \sum_{i \in \N} \sum_{t \in \T} \omega_{i,t} \label{robust:obj}\\
	\text{subject to}
	& \quad \sum_{j \in \N} s_{i,j,t} \leq p_{i,t}  & \forall i \in \N, t \in \T, p_{i,t} \in P_{i,t} \label{robust:con:sent}\\
	& \quad \alpha_{i,t} - b_i \leq \omega_{i,t}  & \forall i \in \N, t \in \T, \alpha_{i,t} \in A_{i,t} \label{robust:con:overflow}\\
	& \quad A_{i,t} = \Big\{ (p_{i,0} - \sum_{t'=1}^{t} d_{i,t'}) + \sum_{j \in \N} s_{i,j,t} \nonumber\\& \quad + \sum_{t'=1}^{t} \big( [1 - \LOS(t - t')] [p_{i,t'} + \sum_{j=1}^{N} (s_{j,i,t'} - s_{i,j,t'})]\big) \nonumber\\& \quad : p_{i,t'} = \bar{p}_{i,t'} + \min \set{0,\xi_{i,t'} \tilde{p}^{-}_{i,t'}} + \max \set{0,\xi_{i,t'} \tilde{p}^{+}_{i,t'}}
	\nonumber\\& \quad
	\forall t' \in \set{1,...,t},
	||\xi_{i,1:t}||_{\infty} \leq 1, ||\xi_{i,1:t}||_{1} \leq \Gamma \Big\} &\forall i \in \N, t \in \T \label{robust:expr:alpha} \\
	& \quad s_{i,j,t} = 0 & \forall (i,j) \in \overline{E(G)}, t \in \T \label{robust:con:graph}\\
	& \quad s_{i,j,t} \geq 0, \quad \omega_{i,t} \geq 0 & \forall i,j \in \N, t \in \T \label{robust:con:0}
\end{align}
\end{subequations}

\iflong
This initial robust formulation is potentially very large relative to the base patient allocation model and is computationally expensive. Assuming that the uncertainty set $P$ is finite, constraints \eqref{robust:con:sent} and \eqref{robust:con:overflow} each represent $|N| \cdot |T| \cdot |P|$ constraints. In comparison, the base model has $|N| \cdot |T|$ constraints, which makes the robust model possibly much more time consuming to solve. If $P$ is instead infinite, then the robust model will have infinite number of constraints. However, it is possible to reformulate the robust model in such a way that it contains the same number of constraints and variables as the base model.

In all the models developed in this work, the number of admitted patients is a constant term in constraints rather than a coefficient on some variable. This means that the uncertainty only exists on the right-hand-side of the linear programming formulation when it is put in standard form, which makes the robust model simpler and easier to solve than a general robust optimization model.
In this case, to ensure that a constraint remains feasible for all realizations in the uncertainty set, we can replace each uncertain constraint with a certain constraint using the worst-case from the uncertainty set of the parameter to compute the right-hand-side.
If the constraint is of the form $ax \leq b$ where $b$ is in some uncertainty set $B$ then the constraint is feasible for all $b \in B$ if and only if it is feasible for $b^{-} = \min B$. Similarly, if the constraint is of the form $ax \geq b$ where $b$ is in some uncertainty set $B$ then the constraint is feasible for all $b \in B$ if and only if it is feasible for $b^{+} = \max B$.

The worst case for constraint \eqref{robust:con:overflow} in the initial robust formulation is the realization from the uncertainty set that yields the largest number of active patients. Therefore, to solve the model, it must be determined which realization yields the largest number of active patients under the constraint $||\xi_i|| \leq \Gamma$.
For a given node $i$ and day $t$ the active patients is a value from the set $A_{i,t}$ defined in constraint \eqref{robust:expr:alpha}. This expression is maximized, generating the worst case, when $\sum_{t^{\prime}=1}^{t} \xi_{i,t} (1-\mathcal{L}(t-t^{\prime})) \tilde{p}^{+}_{i,t^{\prime}}$ is maximized since all other terms are constant or non-positive.
Assuming that each $(1-\mathcal{L}(t-t^{\prime})) \tilde{p}^{+}_{i,t^{\prime}}$ is unique, $\xi_{i,t}$ will be 1 for the $\Gamma$ days with largest $(1-\mathcal{L}(t-t^{\prime})) \tilde{p}^{+}_{i,t^{\prime}}$ and zero otherwise. This results in a closed-form expression for $\xi_{i,t}$ that produces the worst case right-hand-side for constraint \eqref{robust:con:overflow}. Notably, this allows us to solve the robust model without adding any additional constraints or variables to the base LP model.

However, an important consideration for our robust model is that the worst-case scenario from the uncertainty set will not be consistent between constraints.
We therefore must select a worst-case scenario per-constraint. This means that the robust model is more conservative than suggested by the uncertainty budget, but in practice $\Gamma$ can still be selected such that the model has the desired level of robustness.

The reformulation of the robust model is therefore as follows:
\begin{subequations}
\begin{align}
	 \underset{\omega}{\text{minimize}} & \quad \sum_{i \in \N} \sum_{t \in \T} \omega_{i,t} \label{robust_reform:obj}\\
	\text{subject to}
	& \quad \sum_{j \in \N} s_{i,j,t} \leq p^{-}_{i,t},  & \forall i \in \N, t \in \T \label{robust_reform:con:sent}\\
	& \quad \alpha^{+}_{i,t} - b_i \leq \omega_{i,t},  & \forall i \in \N, t \in \T \label{robust_reform:con:overflow}\\
	& \quad s_{i,j,t} = 0 & \forall (i,j) \in \overline{E(G)}, t \in \T \label{robust_reform:con:graph}\\
	& \quad s_{i,j,t} \geq 0 & \forall i,j \in \N, t \in \T \label{robust_reform:con:s0} \\
	& \quad \omega_{i,t} \geq 0 & \forall i \in \N, t \in \T \label{robust_reform:con:o0} \\
	& \quad p^{-}_{i,t} = \bar{p}_{i,t} - \tilde{p}^{-}_{i,t} & \forall i \in \N, t \in \T \label{robust_reform:con:pminus} \\
	& \quad p^{+}_{i,t} = \bar{p}_{i,t} + \tilde{p}^{+}_{i,t} & \forall i \in \N, t \in \T \label{robust_reform:con:pplus} \\
	& \quad \alpha^{+}_{i,t} = (p_{i,0} - \sum_{t'=1}^{t} d_{i,t'}) \nonumber\\& \quad + \sum_{t'=1}^{t} \{ [1 - \LOS(t - t')] [((1 - \zeta_{i,t,t'}) \bar{p}_{i,t'} + \zeta_{i,t,t'} \tilde{p}^{+}_{i,t}) 
	\nonumber\\& \quad
	+ \sum_{j=1}^{N} (s_{j,i,t'} - s_{i,j,t'})]\} &\forall i \in \N, t \in \T \label{robust_reform:con:alphaplus}
\end{align}

Where:
\begin{align}
    &\quad \zeta_{i,t} = \argmax_{\hat{\zeta}_{i,t}} \sum_{t'=1}^{t} \LOS(t-t') \tilde{p}^{+}_{i,t'} \hat{\zeta}_{i,t,t'} & \forall i \in \N, t \in \T \label{robust_reform:subproblem:a}\\
    \text{subject to}
    &\quad \hat{\zeta}_{i,t} \in \set{0,1}^{t}  & \forall i \in \N, t \in \T \label{robust_reform:subproblem:b}\\
    &\quad ||\hat{\zeta}_{i,t}||_{1} = \min\set{\Gamma,t}  & \forall i \in \N, t \in \T \label{robust_reform:subproblem:c}
\end{align}
\end{subequations}
\fi

\subsection{Solution Evaluation} \label{section:methods:evaluation}

We consider the overflow in patient-days, which is defined as the number of patients minus the number of beds, summed over all node-day pairs such that the number of patients is greater than the number of beds: $\sum_{i \in \N} \sum_{t \in \T} \max \set{0, \alpha_{i,t} - b_{i}}$. The overflow is the amount of extra capacity that a node will have to create to properly care for all its patients, or else it must provide sub-standard care or turn patients away. Each of the models directly minimizes this key metric.
A second metric, overflow reduction, is the decrease in total overflow under the model's transfer scheme as a percentage of the baseline overflow, which is the overflow assuming no transfers were made.
Additionally, the size of individual node-day overflows is considered as large overflows on a given node-day necessitate a large amount of surge capacity to avoid inadequate patient care. This means good solutions will have smaller individual node-day overflows, even if this means the total overflow is spread out over more node-days. We therefore evaluate models on the mean, median, and maximum non-zero overflow.
Overflow can also be computed system-wide rather than for a given node-day pair by aggregating all patients and capacity.
This system-wide overflow is not dependent on patient transfers; instead, it helps to contextualize other metrics. A non-zero system-wide overflow indicates that the system as a whole is overburdened and therefore desirable solutions are harder to find, and it represents a lower bound on the total overflow, although in general this bound is not achievable because we only consider transfers of newly admitted patients.

\iflong
\clearpage
\fi

{\TableSpaced
\begin{longtable}[hbt]{
        >{\centering}p{0.09\textwidth}
        >{\arraybackslash}p{0.88\textwidth}
    }
    \caption{General information on the models}
    \label{tab:modelinfo} \\
    \hline \hline
    \setlength{\tabcolsep}{6pt} \\
    Patient Modeling Assumptions & 
    \begin{itemize}
    \setlength\itemsep{0.25em}
        \item We consider only new patients as potential transfers between nodes.
        \item We assume full knowledge about the number of new patients visiting a hospital each day.
        \item We assume that the only resources limiting patient care are hospital beds, which are fixed in number at each node.
        \item We assume that a fixed proportion of hospital beds are available to COVID patients.
        \item All patients have a length of stay governed by a distribution $\LOS_{g}$.
    \end{itemize} \\
    \hline
    \setlength{\tabcolsep}{6pt} \\
    Nurse Modeling Assumptions & 
    \begin{itemize}
    \setlength\itemsep{0.25em}
        \item Nurses may be moved between hospitals.
        \item We have full knowledge of the initial number of nurses in each region, which is constant except for our reallocation.
        \item In addition to bed availability, nurse availability limits patient care.
        \item Some fixed proportion of nurses are available to treat COVID patients.
    \end{itemize} \\
    \hline
    \setlength{\tabcolsep}{6pt} \\
    Sets & \begin{itemize}
        \setlength\itemsep{0.25em}
        \item $\N$: patient treatment nodes, indexed by $n \in \N$
        \item $\T$: modeling days, indexed by $t \in \T = \set{1,2,3,...,T}$
        \item $\G$: patients groups, indexed by $g \in \G$
        \item $\B$: bed types, indexed by $\beta \in \B$
    \end{itemize} \\
    \hline
    \setlength{\tabcolsep}{6pt} \\
    Data &
    \begin{itemize}
        \setlength\itemsep{0.25em}
        \item $p_{g,i,t}$: number of patients admitted to node $i \in \N$ at time $t \in \T$, specified in group $g \in \G$. (zero for $g$s where $indegree(g) \neq 0$)
        \item $p_{g,i,0}$: number of initial patients in group $g \in \G$ at node $i \in \N$ at time $t = 0$
        \item $d_{g,i,t}$: number of initial patients who were discharged from group $g \in \G$ at node $i \in \N$ at time $t \in \T$
        \item $b_{\beta, i}$: number of beds of type $\beta \in \B$ available for COVID patients at node $i \in \N$
        \item $n_i$: initial number of nurses at node $i \in \N$
        \item $G$: directed graph where $V(G) = \N$ and $(i,j) \in E(G)$ if and only if node $i$ may transfer resources to node $j$
        \item $\ell_{g}$: distribution over length of stay for patients in group $g \in \G$
        \item $\LOS_{g}$: cumulative distribution over length of stay for patients in group $g \in \G$
    \end{itemize} \\
    \hline
    \setlength{\tabcolsep}{6pt} \\
    Parameters &
    \begin{itemize}
    \setlength\itemsep{0.25em}
        \item $R_{\text{thresh}}$: load ratio at which balancing penalization begins
        \item $S_{\text{min}}$: minimum number of patients that can be included in a transfer
        \item $T_{\text{switch}}$: minimum number of days a node must wait between sending and receiving patients
        \item $C_{\text{balance}}$: cost coefficient for load-balancing penalty $\phi_{i,t}$
        \item $C_{\text{smooth}}$: cost coefficient for the smoothing penalty $\delta_{i,j,t}$
        \item $C_{\text{sent}}$: cost coefficient for patient transfers $s_{i,j,t}$
        \item $C_{\text{setup}}$: cost coefficient for the transfer indicator $\rho_{i,j}$
        \item $C_{\text{patient}}$: cost coefficient for patient overflow
        \item $C_{\text{nurse}}$: cost coefficient for nurse overflow
        \item $G_{\text{group}}$: a graph where for all $\forall i,j \in \G, i \sim j$ if and only if patients from group $i \in \G$ are transferred to group $j \in \G$
        \item $f$: function mapping $\G$ to $\B$
        \item $Q_{\beta}$: ratio of nurse-days to patient-days for bed type $\beta \in \B$
    \end{itemize}\\
    \hline
    \setlength{\tabcolsep}{6pt} \\
    Variables & \begin{itemize}
        \setlength\itemsep{0.25em}
        \item $s_{g,i,j,t}$: number of patients of patient group $g \in \G$ sent from node $i \in \N$ to node $j \in \N$ at time $t \in \T$
        \item $o_{\beta,i,t}$: dummy variable for patient overflow in bed type $\beta \in \B$ at node $i \in \N$ at time $t \in \T$
        \item $\delta_{i,j,t}$: dummy variable for the absolute difference in the number of patients sent from node $i \in \N$ to node $j \in \N$ between days $t-1$ and $t \in \T \setminus \set{1}$.
        \item $\phi_{\beta,i,t}$: dummy variable for the amount by which patient load ratio exceeds $R_{\text{load}}$ at node $i \in \N$ at time $t \in \T$
        \item $\rho_{i,j}$: binary dummy variable which is equal to 1 if and only if there is a patient transfer between node $i \in \N$ and node $j \in \N$
        \item $\nu_{1,i,t}, \nu_{2,i,t}$: binary dummy variables used to enforce the minimum number of days a node must wait between sending and receiving patients
        \item $\sigma_{i,j,t}$: variable for nurses sent from region $i \in \N$ to region $j \in \N$ at time $t \in \T$
        \item $\theta_{i,t}$: dummy variable for nurse overflow in region $i \in \N$ at time $t \in \T$
    \end{itemize} \\
    \hline
    \setlength{\tabcolsep}{6pt} \\
    Expressions &
    \begin{itemize}
        \setlength\itemsep{0.25em}
        \item $\alpha_{g,i,t}$: expression for the total number of active patients in group $g \in \G$ at node $i \in \N$ at time $t \in \T$
        \item $\chi_{g,i,t}$: expression for the total number of patients entering group $g \in \G$ at node $i \in \N$ at time $t \in \T$
        \item $\gamma_{g,i,t}$: expression for the total number of patients leaving group $g \in \G$ at node $i \in \N$ at time $t \in \T$
        \item $\eta_{i,t}$: expression for the total number of nurses at node $i \in \N$ at time $t \in \T$
        \item $q_{i,t}$: expression for the total nurse demand at node $i \in \N$ at time $t \in \T$
    \end{itemize} \\
    \hline
\end{longtable}
}

Optimal objective function value is also evaluated, which takes into account the overflow as well as the values of the logistical penalties associated with the solution.
For each node-day we also consider the patient load, defined as the number of active patients divided by the number of beds. This load is therefore a measure of the stress on a node that is normalized by its capacity and thus directly comparable between node-day pairs. The maximum load over time is an important metric as it measures how severe the stress on a node will get at the peak, which is what that node will have to prepare for.
Finally, we evaluate solutions on the size of the patient transfers they make as well. Solutions that transfer fewer patients are more efficient and easier to implement in practice. As with overflow, the size of individual transfers matter, so we consider mean, median, and maximum non-zero transfer size.
While these metrics capture many important aspects of the performance of a solution, there are other operational considerations, such as robustness, that are not fully captured. These secondary solution characteristics can be measured using additional metrics, and be qualitatively evaluated using the plots included in Section \ref{section:results}.

\section{Data}
\label{section:data}

The models developed in Section \ref{section:methods} require the following inputs in order to solve the problem of optimal COVID patient redistribution:
\begin{enumerate}
    \item The number of active COVID patients at each location in the study on the first day of the study period.
    \item The proportion of these initial patients that were discharged at each location on each day.
    \item The number of COVID patients admitted to each location on each day, which can optionally be uncertain within some certain interval.
    \item The non-surge capacity of each location available for COVID patients.
    \item The distribution over the length of stay for COVID patients.
    \item A graph representing which pairs of locations are permitted to transfer patients with each other, which can optionally be directed.
\end{enumerate}

In addition to such data, the group models require the capacity per bed type (ICU or ward) and the COVID patient inputs per patient group. Models involving nurse allocation also require nurse supply and demand, which is computed from the number of nurses at each location, the number of hours a nurse works per week, and the number of nurse-hours that should be devoted to each COVID patient per day.

While hospitals in the United States are mandated to report detailed metrics about their COVID patient load and response to the federal government, this data is aggregated at the state level before it is shared publicly, and even then only a small subset of the collected metrics are reported. Beyond this, there is no standardized publicly-available reporting of hospitalization data associated with COVID. However, some state and local governments have elected to release some of this data, enabling us to consider these regions as more adequate case studies for the  methodology proposed in this work.
Due to the lack of standardization in data reporting and the stringent data requirements of our models, data collection, processing, and cleaning proved to be a critical and substantial task. To aid in future studies regarding hospitalizations associated with COVID, we have open-sourced the data we collected and the code to process it, as well as published a website that compiles relevant data sources.
\ifblind
\else
The list of data sources can be found at \url{https://jhu-covid-optimization.github.io/covid-data/}, while the data processing code can be found at \url{https://github.com/flixpar/covid-resource-allocation}.
\fi

In this work, we present case studies based on historical data. In practice, however, this method would need to be implemented prospectively, using projections of the number of COVID-related hospitalizations. Forecasting COVID-related hospitalizations is beyond the scope of this work, but there exist multiple methods for such forecasting at both the hospital level and the state level \citep{CDC2020Forecasting,PennMedicine2020}. The reader is referred to such models for more information.

\iflong
\begin{table}[ht]
\TableSpaced
\small
    \caption{Data sources used throughout this work.}
    \label{fig:data:sources}
    \begin{tabular}{>{\arraybackslash}p{0.25\textwidth}>{\arraybackslash}p{0.3\textwidth}>{\arraybackslash}p{0.4\textwidth}}
    \hline\hline
        Description & Source & Source URL \\
        \hline 
         Florida Hospitalizations & Florida AHCA & \tiny{\url{https://ahca.myflorida.com/COVID-19_alerts.shtml}}\\
        New Jersey Hospitalizations & New Jersey DOH & \tiny{\url{https://services7t:.arcgis.com/Z0rixLlManVefxqY/arcgis/rest/services/PPE_Capacity/FeatureServer/0}} \\
        Texas Hospitalizations & Texas DOH & \tiny{\url{https://dshs.texas.gov/coronavirus/additionaldata/}} \\
        Hospital Beds & Definitive Healthcare &  \tiny{\url{https://coronavirus-resources.esri.com/datasets/1044bb19da8d4dbfb6a96eb1b4ebf629_0}} \\
        COVID-19 LOS Distribution & \cite{lewnard2020incidence} & \tiny{\url{https://www.bmj.com/content/369/bmj.m1923}}
        \\\hline
    \end{tabular}
\end{table}

\subsection{COVID-Related Hospitalizations}
The primary input data for each of the models developed in this work includes the number of initial active COVID patients, the proportion of initial patients discharged each day, and the number of newly admitted COVID patients each day, for each location. While such data may be known or projected for a specific case study, usually, the number of active COVID patients is more readily available. Therefore, we estimate the number of admitted patients and the proportion of initial patients discharged on each day from the reported number of active COVID patients. Specifically, for each location, given a candidate solution for the number of admitted patients and the proportion of initial patients discharged on each day, the number of active patients on each day is computed using constraint \eqref{base:expr:alpha} assuming zero patients are transferred. The objective function is then the $\ell_2$-distance between the computed number of active patients on each day and the reported number of active patients, which is approximately minimized over the candidate admitted and discharged solutions using random search optimization.

For the specific case studies considered in this work, different sources of data are considered to obtain COVID-related hospitalization data (see Figure \ref{fig:data:sources}).
The first case study we consider is New Jersey. We obtained the number of active COVID patients from the New Jersey Department of Health. This dataset also reports the number of COVID patient admissions, discharges, and bed availability by hospital. However, the reported number of admissions and discharges are not entirely consistent with the reported number of active COVID patients, and therefore we decided to estimate these quantities from the reported active COVID patients using the above method.
The second case study we consider is Florida. The Florida Agency for Healthcare Administration reports the number of active patients by hospital, but only reports the number of active COVID patients by county. We therefore estimate the number of active COVID patients at each hospital by multiplying the total number of patients at that hospital by the proportion of all patients that have COVID in the county where the hospital is located.
We have also found that the number of active patients in Florida has sporadic single-day outliers which appear to be reporting errors. To resolve this problem, we run a simple outlier detection and correction algorithm. For each day and location we compute the median and median absolute deviation (MAD) of active patients over the surrounding five-day period. If the number of active patients deviates from the median by more than ten times the MAD then we consider it to be an outlier, and replace it with the median.
The final case study we consider is Texas. The Texas Department of Health reports the number of active COVID patients by Trauma Service Area (TSA). This data does not appear to require correction, but it does not include patient admissions, so we estimate the COVID patient admissions and discharges using the method above.

\subsection{Hospital Capacity} \label{section:data:beds}
To determine the capacity of each location in our case studies to care for COVID patients we used the Definitive Healthcare USA Hospital Beds dataset (see Table \ref{fig:data:sources}) to get the number of ICU and ward beds at each hospital in the United States. This data was collected before the onset of COVID in the United States, so it represents the ordinary hospital capacity rather than the surge capacity. We assumed that 35\% of ward beds and 50\% of ICU beds could be made available for COVID patients. These parameters are estimated from the state-level COVID hospitalization data reported by the Center for Disease Control (CDC). In the patient allocation models we assume that the number of beds that can be made available for COVID patients determines the capacity. In the nurse allocation model we instead assume that the number of nurses is the limiting factor for capacity instead.

\subsection{Patient Length of Stay} \label{section:data:los}
The distribution over the length of stay (LOS) in ward and ICU beds for COVID patients is estimated by \cite{lewnard2020incidence}. We use a Weibull($\lambda=12.88, k=1.38$) distribution for ward patients and a Weibull($\lambda=13.32, k=1.58$) for ICU patients. See Figure \ref{fig:data:los}. We discretize these distributions for use in the model. In the group model, it is required that patients stay for two days in a non-ICU bed before ICU admittance \citep{wunsch2011comparison} and five days after ICU discharge. \citep{tiruvoipati2017intensive}

\begin{figure}
    \centering
    \includegraphics[width=0.65\textwidth]{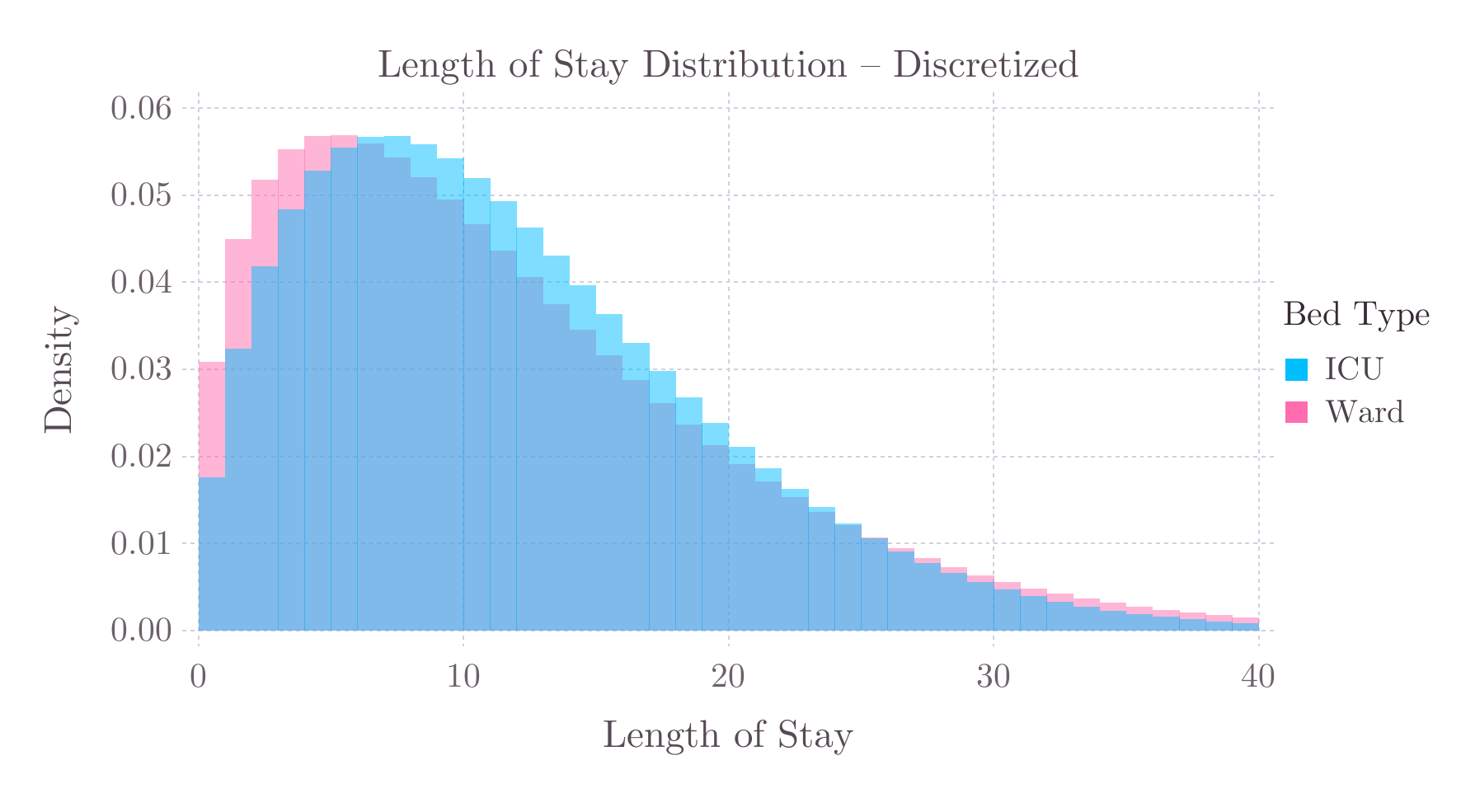}
    \caption{Distribution over the length of stay (LOS) for COVID patients by bed type, discretized by day.}
    \label{fig:data:los}
\end{figure}

\subsection{Adjacency Graph} \label{section:data:adj}
Our models constrain the solution so that patients are transferred only along the edges of a graph, which is an input to the model. In practice, we use this constraint to set an upper limit on the distance that a patient can be transferred by adding and edge between two locations in the graph if they are within a certain distance of each other. Given the latitude and longitude of the centroid of each location, we compute the distance matrix using the Haversine formula. We then threshold this distance matrix to determine the binary adjacency matrix for the patient transfer graph.

\subsection{Nurse Availability} \label{section:data:nurses}

In our New Jersey case study we validate our nurse and combined nurse plus patient models, which require data on the number of nurses that work at each hospital as well as the nurse demand.
We estimate the number of nurses at each hospital from the New Jersey Department of Health's Hospital Patient Care Staffing Report, which reports the mean number of patients per nurse.
Multiplying this quantity by the mean number of patients-days per week we get the mean number of nurses-days per week. Assuming nurses work 36 hours per week on average, we multiply by 24/36 to get the number of total nurses.
We estimate nurse demand from the number of active patients. The number of active patients times the number of nurses per patient yields the number of nurse-days required. We then compute the number of nurses required to supply that many nurse days per actual day, assuming that nurses work 36 hours per week on average.

In general, the number of nurses at each hospital would be known by the decision-maker using this model, and therefore it would not have to be estimated this way.

\subsection{Data Collection Methods and Resources} \label{section:data:collection}
In the process of collecting data for this work we found a great deal of data pertaining to COVID or likely relevant to studies on COVID, from COVID forecasts to nursing workforce surveys to health and hospital statistics.
While there is a lot of available data for studies on COVID, there are also many gaps in the available data, which makes it difficult to know what data can be found and used without an extensive search. Much of the data is also not easily locatable as many sources show only a dashboard summarizing the available data, but the complete data can be accessed as well.
We therefore have collected a list of over fifty data sources related to COVID that we have used in this work or think are likely to be useful to other researchers.
We also collect basic metadata about each data source, including a category, a description, an indication of whether the data is geographical, and the unit resolution. We are making this list public on our website so that other researchers may use it.
\ifblind
\else
It can be found at: \url{https://jhu-covid-optimization.github.io/covid-data/}.
\fi
We are also releasing the code used for this study, which contains scripts used to download, process, clean, and load the data we worked with.
\ifblind
The URL for the website as well as the code is omitted to maintain the anonymity of the authors during review.
\else
This code can be found at: \url{https://github.com/flixpar/covid-resource-allocation}.
\fi

\fi

\section{Results}
\label{section:results}
In this section, we apply the models developed in Section \ref{section:methods} to a number of case studies in order to showcase the models' applicability to real examples, validate them, and demonstrate their effectiveness. In particular, we study the peak of the first wave of the COVID pandemic in New Jersey and Texas. Miami-Dade county is also studied in the online supplement to this paper. We first evaluate each of the models we have developed to justify the inclusion of the additions we made to the base model. We then apply our final model to each of our case studies to demonstrate the potential effectiveness and impact of our approach.

\subsection{Model Validation}

To validate the models developed in this work, we first apply them to our New Jersey case study. New Jersey was selected for this task because the New Jersey Department of Health has published detailed data on COVID patients and bed availability at the hospital level.
Specifically, we consider each of the 75 hospitals in the state of New Jersey from April $5^\text{th}$ to June $15^\text{th}$, 2020.
In Table \ref{fig:results:valid:metrics}, we evaluate each of the non-group patient allocation models according to the metrics defined in Section \ref{section:methods:evaluation}. We also include a "no transfer" model that does not transfer patients as a baseline against which to compare the effectiveness of our models.
We note that without transfers, there is a significant overflow of COVID patients in patient-days, yet, as can be seen in Figure \ref{fig:results:nj:active_total}, the number of total active COVID patients is always less than the capacity of the whole system. Because of this we should expect to see a large reduction in the total overflow due to our models.
Compared to the baseline, all models result in an overflow reduction of at least 28,862 patient-days, or 86.4\%, which is a significant improvement for the hospital system.

\iflong
\begin{figure}[tb]
    \centering
    \includegraphics[width=0.55\textwidth]{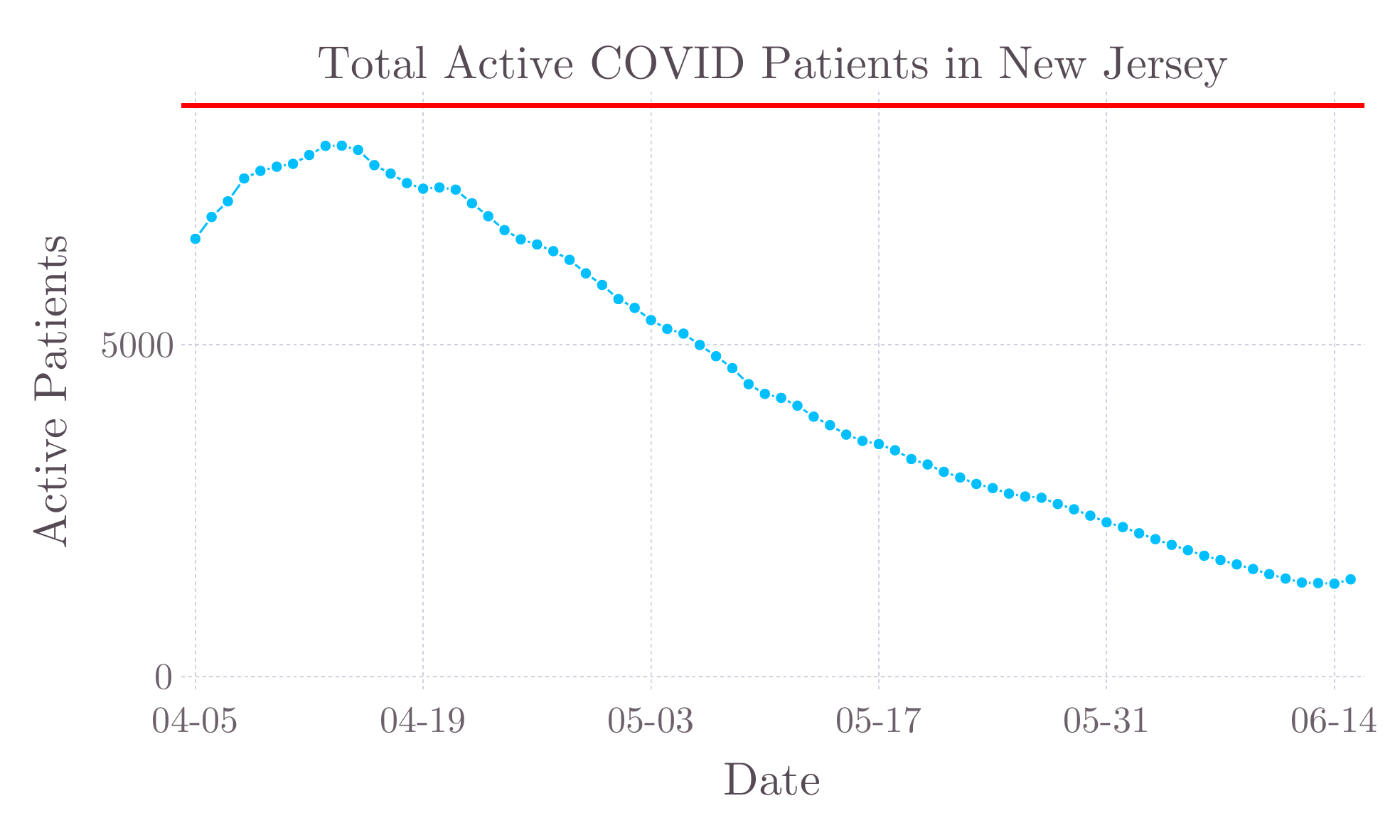}
    \caption{Total active COVID patients by day in New Jersey.}
    \label{fig:results:valid:active_total}
    \label{fig:results:nj:active_total}
\end{figure}
\fi

\begin{figure}[ht]
    \centering
    \includegraphics[width=1.0\textwidth]{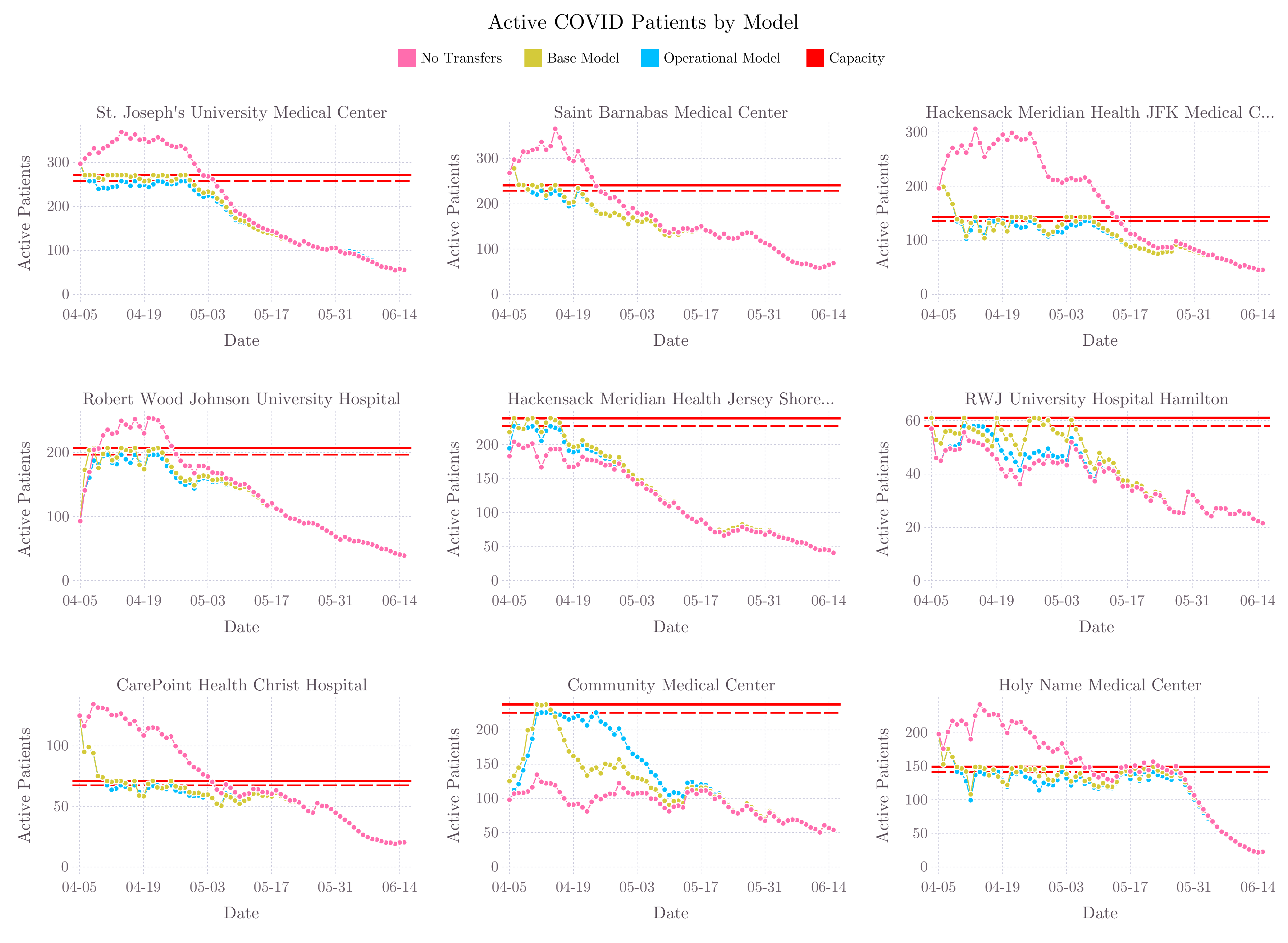}
    \caption{Active patients by model for a representative sample of nine hospitals in New Jersey.}
    \label{fig:results:valid:active_subset}
    \vspace{-1.5em}
\end{figure}

The base model has the least restrictive constraints and penalties, and therefore is able to achieve the largest decrease in overflow, which is the primary goal. It also performs well according to other metrics including mean and median non-zero overflow. However, it does have some undesirable solution characteristics as compared to the other models, such as a higher percentage of hospital-days with an overflow, more total patients transferred, and a larger maximum-sized transfer. Figure \ref{fig:results:valid:active_subset} demonstrates that the base model often brings hospitals that were well under capacity all the way up to capacity as well, which in practice they will likely be unwilling to accommodate.
The operational model addresses many of these issues, at the cost of greater total overflow. In particular, it transfers fewer patients overall and the mean, median, and maximum transfer size are all significantly smaller than those of the base model, which would make the transfers more operationally feasible.

\begin{table}[ht]
    \TableSpaced
    \centering
    \caption{Evaluation of the performance of each of the non-group patient allocation models in our New Jersey case study.}
    \label{fig:results:valid:metrics}
    \begin{tabular}{
    >{\arraybackslash}p{0.3\textwidth}
    |
    >{\raggedleft\arraybackslash}p{0.10\textwidth}
    >{\raggedleft\arraybackslash}p{0.10\textwidth}
    >{\raggedleft\arraybackslash}p{0.11\textwidth}
    >{\raggedleft\arraybackslash}p{0.13\textwidth}
    >{\raggedleft\arraybackslash}p{0.13\textwidth}
}
\hline\hline
  & No & Base & Operational & Base & Operational\\
  & Transfers & Model & Model & Robust Model & Robust Model\\
\hline
Overflow & $33406$ & $3800$ & $3930$ & $4173$ & $4544$\\
Overflow Reduction & $0.0\%$ & $88.62\%$ & $88.24\%$ & $87.51\%$ & $86.4\%$\\
Median Non-Zero Overflow & $27.5$ & $8$ & $11$ & $21$ & $9$\\
Mean Non-Zero Overflow & $44.7$ & $27.7$ & $30.2$ & $35.4$ & $26.6$\\
Max Non-Zero Overflow & $244$ & $220$ & $220$ & $220$ & $220$\\
Median Load & $40.98\%$ & $51.8\%$ & $52.0\%$ & $53.0\%$ & $53.52\%$\\
Mean Load & $62.15\%$ & $55.2\%$ & $54.91\%$ & $55.18\%$ & $56.54\%$\\
Max Load & $2000.0\%$ & $400.0\%$ & $400.0\%$ & $400.0\%$ & $900.0\%$\\
Percent Of Hospital-Days With An Overflow & $13.85\%$ & $2.54\%$ & $2.41\%$ & $2.19\%$ & $3.17\%$\\
\hline
Total Patients Transferred & $0$ & $4298$ & $3932$ & $4894$ & $4255$\\
Percent Of Patients Transferred & $0.0\%$ & $14.25\%$ & $13.04\%$ & $16.23\%$ & $14.11\%$\\
Median Non-Zero Transfer & $0$ & $7$ & $1$ & $5$ & $2$\\
Mean Non-Zero Transfer & $0$ & $9.7$ & $2.4$ & $7.8$ & $2.6$\\
Max Non-Zero Transfer & $0$ & $55$ & $19$ & $52$ & $25$\\
Percent Of Hospital-Days With A Transfer & $0.0\%$ & $12.81\%$ & $24.31\%$ & $17.35\%$ & $24.65\%$\\
\hline
\end{tabular}
    \vspace{-1.5em}
\end{table}

\begin{table}[ht]
    \TableSpaced
    \caption{Evaluation of the performance of each of the group patient allocation models in our New Jersey case study.}
    \label{fig:results:valid:metrics_groupmodel}
    \centering
    \begin{tabular}{
    >{\arraybackslash}p{0.3\textwidth}
    |
    >{\raggedleft\arraybackslash}p{0.09\textwidth}
    >{\raggedleft\arraybackslash}p{0.09\textwidth}
    >{\raggedleft\arraybackslash}p{0.11\textwidth}
    |
    >{\raggedleft\arraybackslash}p{0.09\textwidth}
    >{\raggedleft\arraybackslash}p{0.09\textwidth}
    >{\raggedleft\arraybackslash}p{0.11\textwidth}
}
\hline\hline
                                          &                         & ICU        &                   &                         & Ward       &                    \\
\hline
                                          & No Transfers      & Base Group Model & Operational Group Model & No Transfers & Base Group Model & Operational Group Model  \\
\hline
Overflow                                  & $13452$                 & $1129$     & $3189$            & $37068$                 & $18270$    & $28912$            \\
Overflow Reduction                        & $0.0\%$                 & $91.6\%$   & $76.3\%$          & $0.0\%$                 & $50.7\%$   & $22.0\%$           \\
Median Non-Zero Overflow                  & $8$                     & $1$        & $2$               & $25$                    & $16$       & $26$               \\
Mean Non-Zero Overflow                    & $12$                    & $3.6$      & $4.1$             & $37.3$                  & $24$       & $37.3$             \\
Max Non-Zero Overflow                     & $89$                    & $37$       & $37$              & $231$                   & $148$      & $226$              \\
Median Load                               & $50.0\%$                & $78.6\%$   & $75.0\%$          & $62.7\%$                & $73.3\%$   & $68.7\%$           \\
Mean Load                                 & $84.4\%$                & $71.9\%$   & $73.2\%$          & $82.0\%$                & $76.6\%$   & $79.5\%$           \\
Max Load                                  & $735.7\%$               & $433.3\%$  & $433.3\%$         & $591.5\%$               & $414.9\%$  & $580.9\%$          \\
Percent of 
Hospital-Days with an Overflow & $32.5\%$                & $9.1\%$    & $22.6\%$          & $28.8\%$                & $22.0\%$   & $22.5\%$           \\
\hline
Total Patients Transferred                & $0$                     & $1432$     & $1016$            & $0$                     & $3060$     & $1395$             \\
Percent of Patients Transferred           & $0.0\%$                 & $39.4\%$   & $27.9\%$          & $0.0\%$                 & $17.6\%$   & $8.0\%$            \\
Median Non-Zero Transfer                  & $0$                     & $1$        & $1$               & $0$                     & $3$        & $1$                \\
Mean Non-Zero Transfer                    & $0$                     & $2$        & $1.4$             & $0$                     & $4.8$      & $2.2$              \\
Max Non-Zero Transfer                     & $0$                     & $13$       & $13$              & $0$                     & $52$       & $36$               \\
Percent of Hospital-Days with a Transfer  & $0.0\%$                 & $29.9\%$   & $25.2\%$          & $0.0\%$                 & $24.3\%$   & $18.9\%$           \\
\hline
\end{tabular}
    \vspace{-1.5em}
\end{table}
Parameter and constraint selection is important for the operational model and their optimal choices are highly dependent on input data and model use case.
In general, we believe that solution characteristics such as fewer transferred patients, no large spikes in the number of transferred patients, and maintaining a cushion between the number of active patients and the capacity are important because they make the solution more feasible in practice.
Adding small penalties on the total amount of patients transferred and the smoothness in the number of patients transferred between days also seems to improve the solution quality while not significantly increasing the running time of the solver.
Therefore, $C_{\text{sent}} = 0.01,$ $C_{\text{smooth}} = 0.01$ are used in the operational model.
We also include \eqref{opt:con:shortage1} and \eqref{opt:con:shortage2}, which ensure that hospitals are not sent over or further over capacity by our model solution, as we believe this solution characteristic will be important to decision-makers.
On the other hand, setting a lower bound on the size of a patient transfer or adding a setup cost between hospitals may only make sense in some situations, and requires adding integer variables to the model which drastically increases the time it takes to solve.
In our operational model we employ only the optional constraints and penalties that have broad applicability, however, the other optional constraints and penalties remain in our model to provide flexibility to practitioners.

In Figure \ref{fig:results:valid:active_subset} we compare the number of active patients at a subset of the hospitals under the base model, the operational model, and without transfers. It can be seen in this figure how the models perform in utilizing excess capacity in the system to accommodate the surge in demand on regional nodes. In this example, both the models achieve their goal of eliminating the overflow in all the hospitals after a couple of days of redistributing the patients. The base model on occasion sends hospitals that were not experiencing an overflow over capacity, which in practice hospitals may not be willing to accommodate. The operational model on the other hand maintains a 5\% buffer between the number of active patients and the capacity for each hospital-day.

\begin{figure}[ht]
    \centering
    \includegraphics[width=0.9\textwidth]{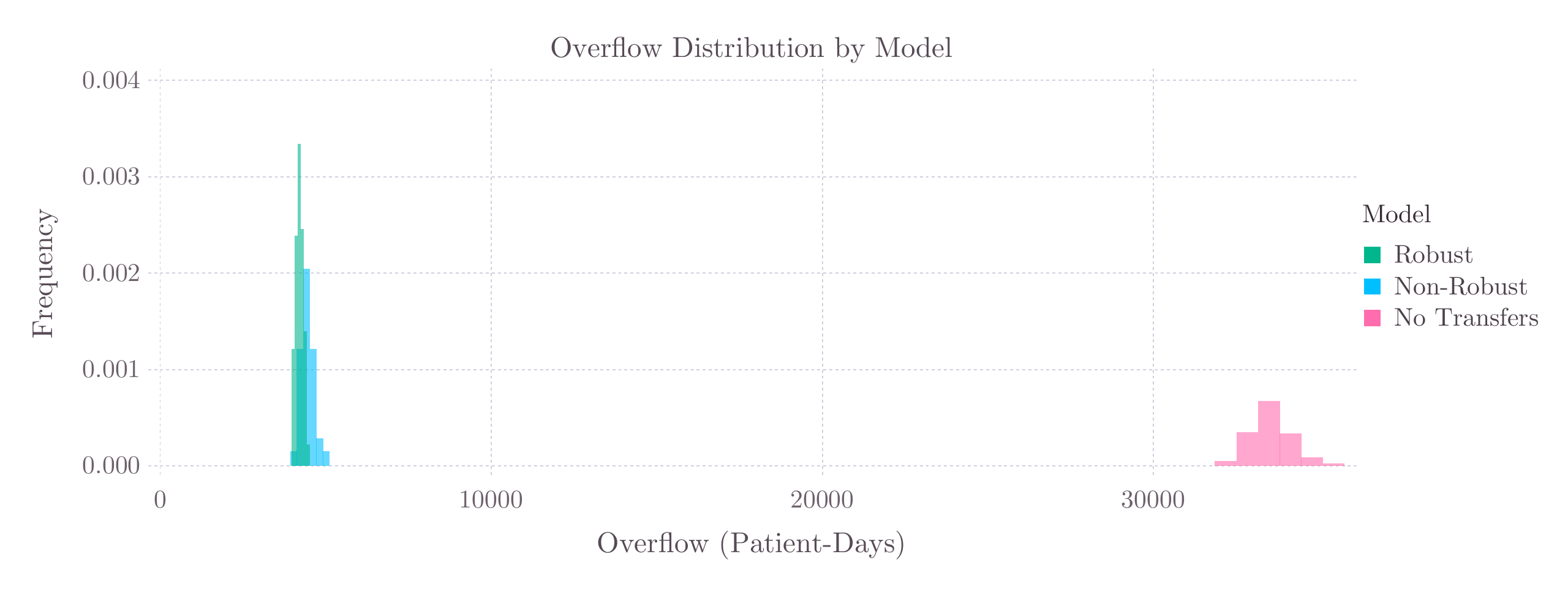}
    \caption{Distribution over the total overflow for all hospitals in New Jersey for three models under scenarios sampled uniformly from the uncertainty set.}
    \label{fig:results:valid:robust_dist}
    \vspace{-1.5em}
\end{figure}

\begin{figure}[ht]
    \centering
    \includegraphics[width=0.9\textwidth]{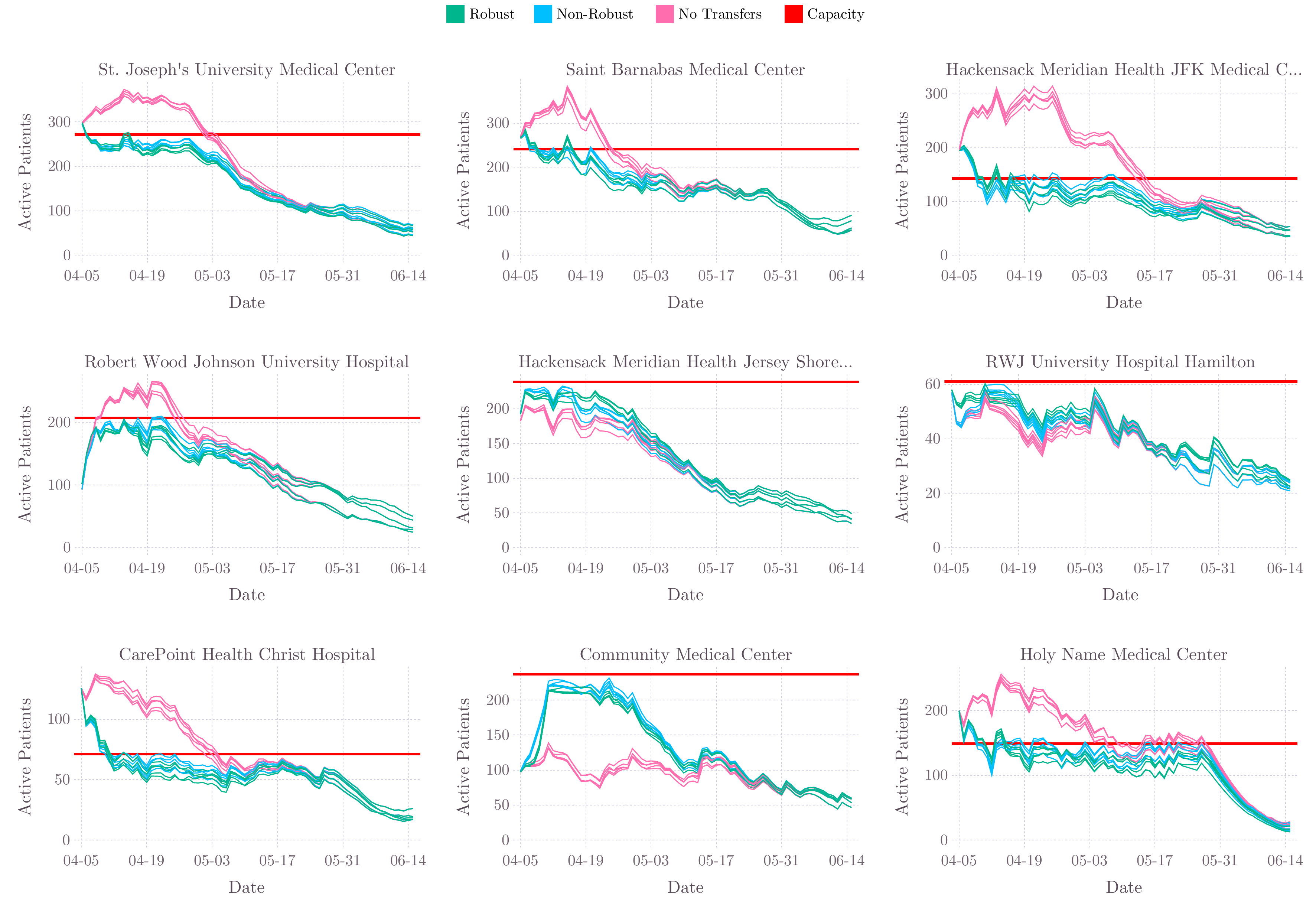}
    \caption{Active COVID patients over time under five scenarios sampled uniformly from the robust uncertainty set at a representative sample of nine hospitals in New Jersey.}
    \label{fig:results:valid:active_robust}
    \vspace{-1.5em}
\end{figure}

The robust model also addresses some of the shortcomings of the base model, however, its primary purpose is to protect the feasibility of the solution against different prospective.
In Figure \ref{fig:results:valid:active_robust} we see the results of the base, robust, and no-transfer models number with active patient counts sampled uniformly from the uncertainty set. In this figure we can see scenarios where the base model sends a hospital over capacity while the robust model does not because the base model only considers the nominal number of admitted patients. The safety associated with the robust model will be its main value to decision-makers, who must account for uncertainty and minimize the negative impact of their actions and policies. It is also important to note that all models perform at least as well as the no-transfer model in all possible scenarios within the uncertainty set. Figure \ref{fig:results:valid:robust_dist} plots the outcome in terms of total overflow from 400 scenarios sampled uniformly from the uncertainty set, and demonstrates that for all of these cases both the base model and the robust model perform better than the baseline.
In this figure we also notice that the robust model consistently performs worse than the base model, which is an unavoidable result of robustness shrinking the the LP's feasible region. Fortunately, we see that this gap is relatively small compared with the gap between each model and the no-transfer model.

The group patient allocation models add the capability to differentiate between patients in different care-paths and to have multiple bed types. We ran the group patient allocation models on the New Jersey hospital-level data at a higher resolution, considering both ward and ICU patients, and in Figure \ref{fig:results:valid:metrics_groupmodel} we compare the results for each group.
These figures show that the base group model is able to achieve a smaller overflow yet has to transfer more patients in both the ICU and ward blocks than either the no-transfer model or the group operational model.
Interestingly, we see that it is not possible to reduce overflow as much in the ward as in the ICU. This is because all patients that visit the ICU also stay in the ward for seven days (two before ICU, five after ICU) but can only be transferred at the beginning of their stay, which means that the number of patient-days in the ward is much higher and there is less flexibility to transfer those extra patient-days.

The last set of models we analyze here are the nurse allocation models. There is a large shortage of nurses-days in the case study we consider because of the increased patient load. Since nurses are critical to properly care for all patients, represents a huge issue. A nurse shortage can be alleviated by increasing the number of hours that each nurse works, utilizing part-time nurses more, or decreasing the care given to patients, however, these solutions are clearly undesirable because they place a larger burden on nurses and have the potential to hurt patient outcomes. Transferring nurses is an appealing alternative to these choices. We see in Table \ref{fig:results:valid:metrics_nurses} that our nurse allocation model is able to alleviate much of the nurse shortage.

\begin{figure}[ht]
    \centering
    \includegraphics[width=0.75\textwidth]{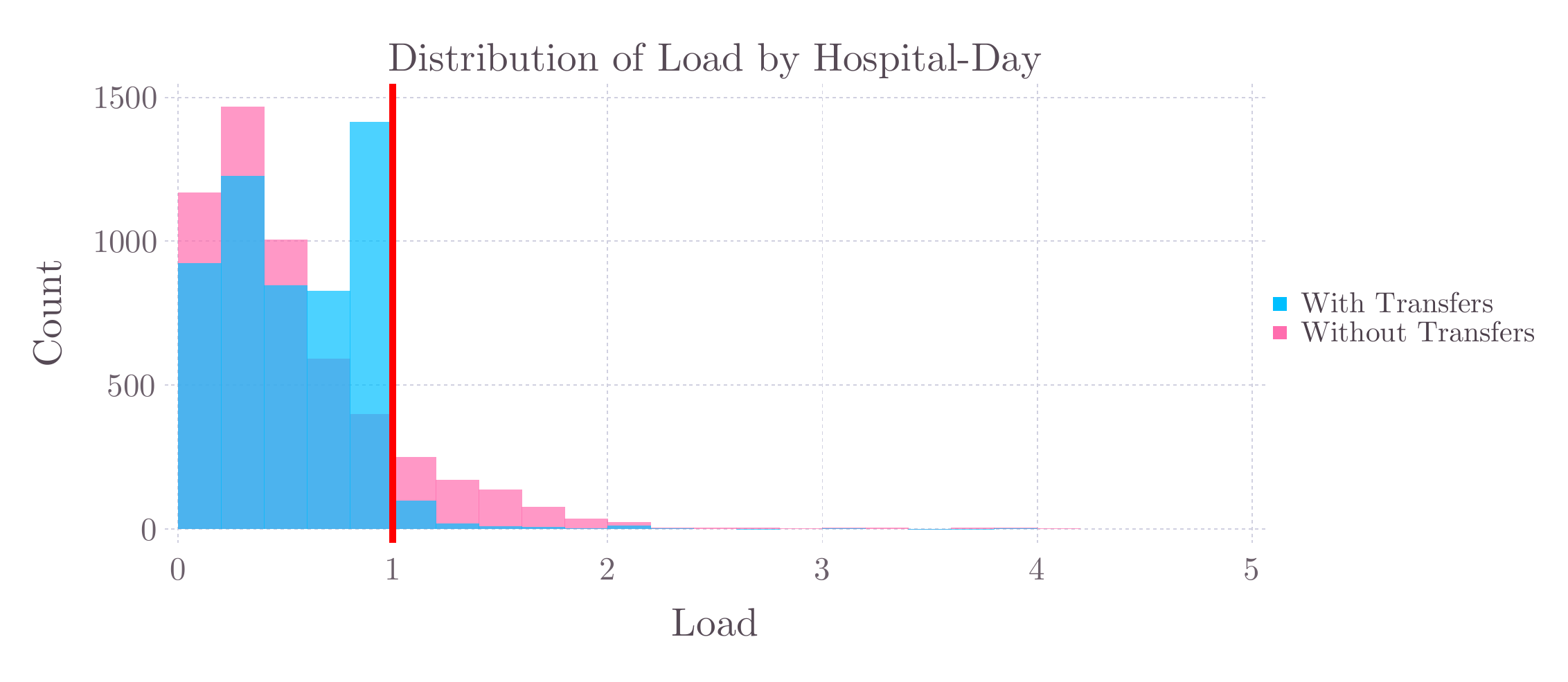}
    \caption{Distribution over the COVID patient load by hospital-day in New Jersey with the operational model.}
    \label{fig:results:valid:overflow_dist}
    \vspace{-1.5em}
\end{figure}

\begin{table}[ht]
    \TableSpaced
    \centering
    \caption{Evaluation of the performance of each of the nurse allocation model in our New Jersey case study.}
    \label{fig:results:valid:metrics_nurses}
    \begin{tabular}{l|rrr}
\hline \hline
  & No Transfers & Base Nurses Model & Operational Nurse Model\\
\hline
Shortage & $55119$ & $34354$ & $34629$\\
Shortage Reduction & $0.0\%$ & $37.7\%$ & $37.2\%$\\
Median Non-Zero Shortage & $38$ & $22$ & $22$\\
Mean Non-Zero Shortage & $47$ & $32$ & $31.9$\\
Max Non-Zero Shortage & $145$ & $202$ & $203$\\
Percent Of Hospital-Days With A Shortage & $67.9\%$ & $62.0\%$ & $62.8\%$\\
\hline
Total Nurse Transfers & $0$ & $8851$ & $1025$\\
Median Non-Zero Transfer & $0$ & $8$ & $1$\\
Mean Non-Zero Transfer & $0$ & $18.2$ & $3.1$\\
Max Non-Zero Transfer & $0$ & $195$ & $74$\\
Percent Of Hospital-Days With A Transfer & $0.0\%$ & $32.6\%$ & $22.4\%$\\
\end{tabular}
    \vspace{-1.5em}
\end{table}

\iflong
\begin{figure}[ht]
    \centering
    \includegraphics[width=0.9\textwidth]{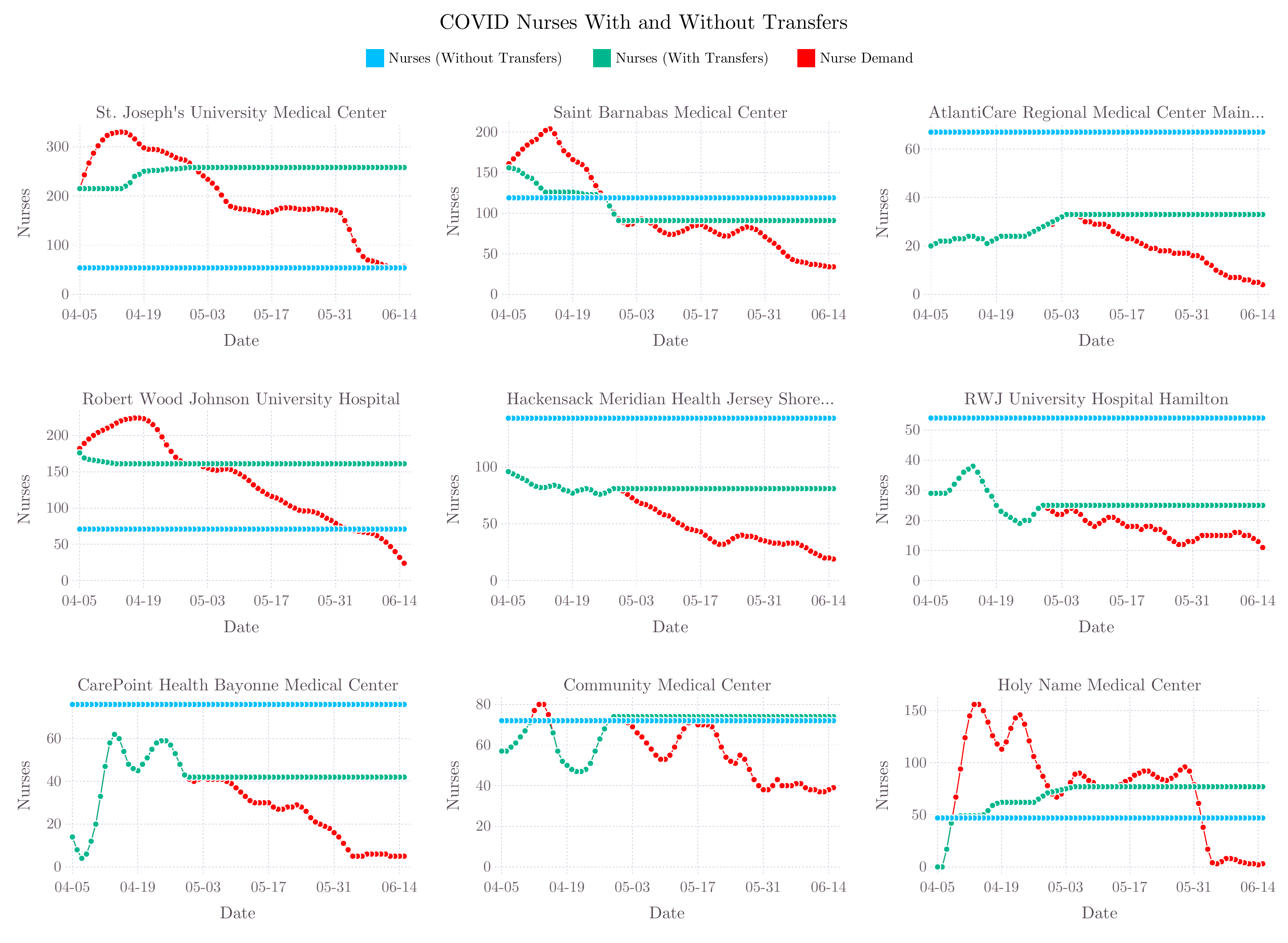}
    \caption{Nurse supply and demand for COVID patients at a representative subset of nine hospitals in New Jersey.}
    \label{fig:results:valid:nurse_active}
\end{figure}
\fi

\iflong
\subsubsection{Model Validation: Load Balancing}
To compare the results of the overflow minimization objective to the load balancing objective we analyze the New York City Health + Hospitals (NYC H+H) System over the period March $27^\text{th}$ to June $1^\text{st}$, 2020.
The NYC H+H System is the largest public healthcare system in the United States, and operates 11 hospitals in New York City. During the peak of the pandemic in early April the system faced an extreme surge in patients, causing all 11 hospitals to go over capacity.
In Figure \ref{fig:results:nychh:active_total} we see that the system as a whole significantly exceeded capacity for an entire month, which makes it a prime candidate for load balancing.
In fact, based on Table \ref{fig:results:nychh:metrics}, the base model with the overflow minimization objective was able to reduce overflow by just 6.8\% because there were few days when some hospitals in the system had excess capacity while others did not. While even a small reduction in the amount of surge capacity that must be created could be useful, load balancing is likely a more useful objective in this context because it means that none of the hospitals has to face an extreme surge while others do not.

Without transfers, the maximum load at a hospital is 205\%, while with load balancing this shrinks to 160\%.
However, there can be a downside to balancing the load as evenly as possible, which is that it often involves far more patient transfers. In general, balancing load as evenly as possible involves micro-adjustments and means patients must be transferred even when all hospitals are over or under capacity, which greatly increases the number of transfers to be made. In Table \ref{fig:results:nychh:metrics} we observe that the base model with a load balancing objective gives a solution that involves the transfer of 1,065 patients whereas the overflow minimization objective leads to just 394 patient transfers. At the same time, the mean size of a transfer decreases to 2.9 patients from 6 patients when load balancing.
Some of this issue is mitigated when adding operational penalties and constraints to the model. The operational load balancing model transfers 453 patients, which is far more similar to solution of the overflow minimization model. The trade-off is that maximum load increases to 174\%.

Evidently, the load-balancing objective can be useful in some situations, although there are trade-offs to using it over the overflow minimization objective.

\begin{figure}
    \centering
    \includegraphics[width=0.5\textwidth]{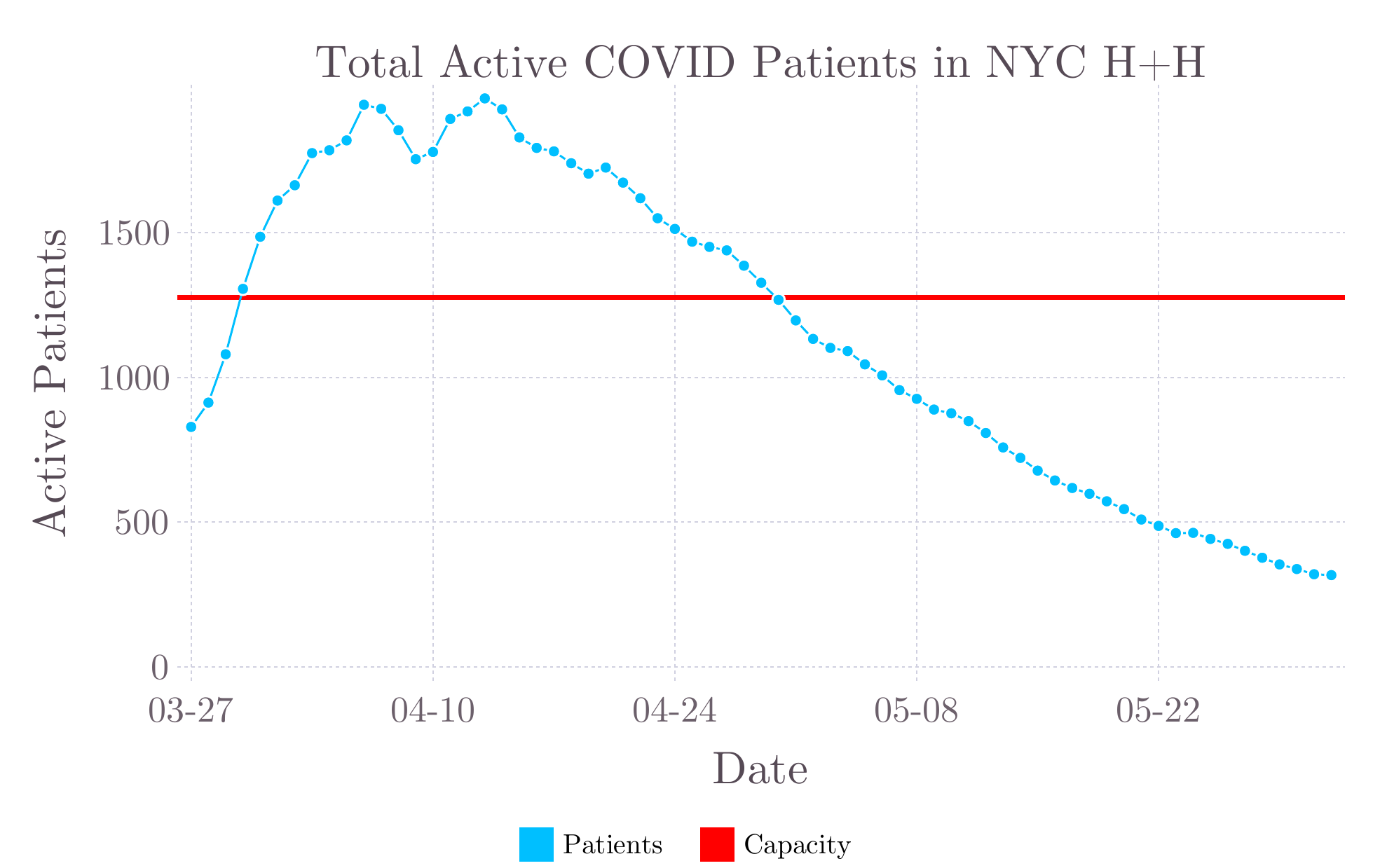}
    \caption{Total active COVID patients (blue) versus the bed capacity for COVID patients (red) in the New York City Health + Hospitals System.}
    \label{fig:results:nychh:active_total}
\end{figure}

\begin{table}
    \TableSpaced
    \centering
    \caption{Evaluation of the performance of each model in the NYC H+H case study.}
    \label{fig:results:nychh:metrics}
    \begin{tabular}{
    >{\arraybackslash}p{0.30\textwidth}
    |
    >{\raggedleft\arraybackslash}p{0.15\textwidth}
    >{\raggedleft\arraybackslash}p{0.15\textwidth}
    >{\raggedleft\arraybackslash}p{0.15\textwidth}
    >{\raggedleft\arraybackslash}p{0.15\textwidth}
}
\hline\hline
 & No Transfers & Overflow Minimization (Base) & Load Balancing (Base) & Load Balancing (Operational) \\
\hline
Overflow & $13954$ & $13011$ & $13560$ & $13275$\\
Overflow Reduction & $0.0\%$ & $6.8\%$ & $2.8\%$ & $4.9\%$\\
Median Load & $88.1\%$ & $95.9\%$ & $91.8\%$ & $93.6\%$\\
Mean Load & $91.2\%$ & $91.6\%$ & $93.1\%$ & $92.5\%$\\
Maximum Load & $205.3\%$ & $190.9\%$ & $160.1\%$ & $173.9\%$\\
\hline
Total Patients Transferred & $0$ & $394$ & $1065$ & $453$\\
Percent Of Patients Transferred & $0.0\%$ & $5.7\%$ & $15.5\%$ & $6.6\%$\\
Percent Of Hospital-Days With A Transfer & $0.0\%$ & $12.9\%$ & $61.3\%$ & $35.7\%$\\
Median Non-Zero Transfer & $0$ & $3$ & $2$ & $2$\\
Mean Non-Zero Transfer & $0$ & $6$ & $2.9$ & $2.2$\\
Max Non-Zero Transfer & $0$ & $29$ & $26$ & $8$\\
\end{tabular}
\end{table}

\begin{figure}
    \centering
    \includegraphics[width=1.0\textwidth]{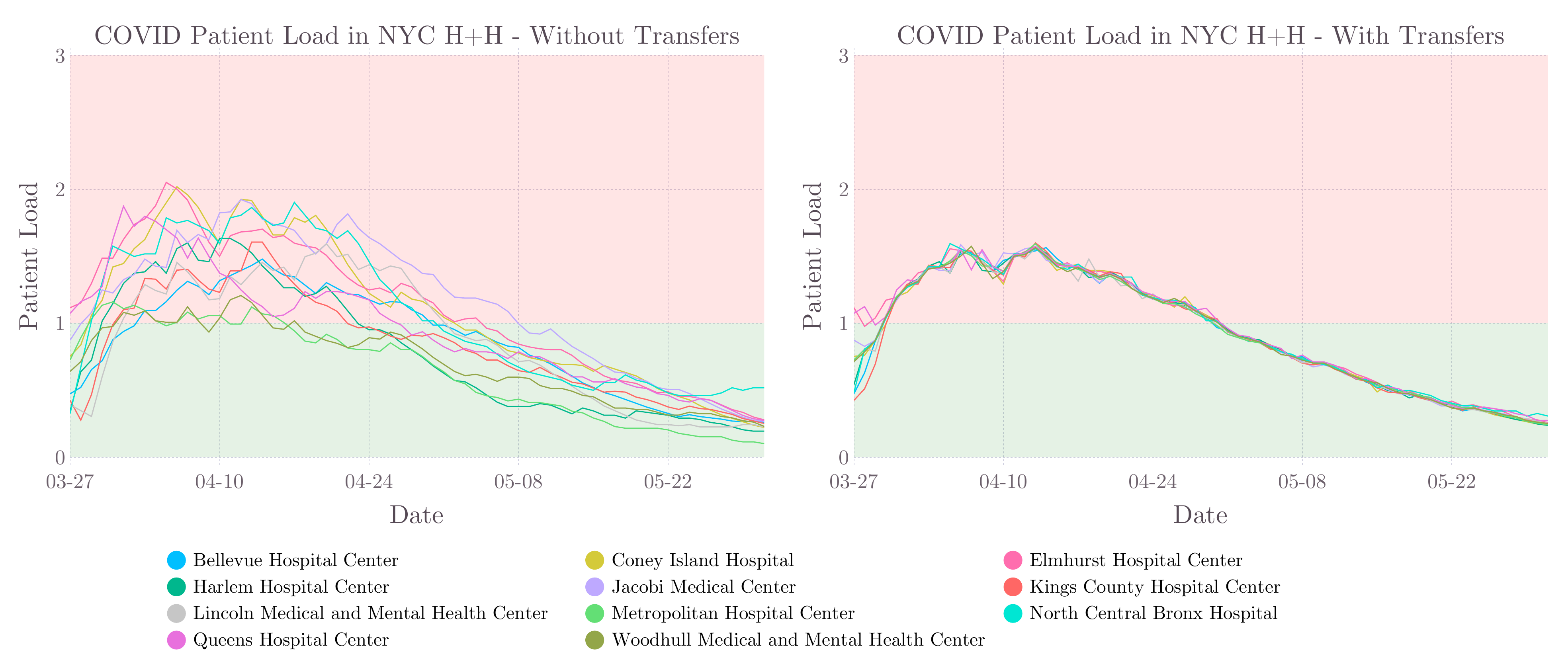}
    \caption{Normalized load over time in the NYC H+H case study using the base demand redistribution model with the load-balancing objective.}
    \label{fig:results:nychh:loadbalance_base:load}
\end{figure}
\fi

\subsection{Patient Redistribution Case Studies}
So far, we have validated our models on a case study in New Jersey. This section focuses on the solution and operational insights; we apply the operational model to three primary case studies, one in New Jersey, one in Texas, and one in Florida's Miami-Dade County, to investigate the solution and its consequences.

\iflong
\begin{table}[ht]
    \TableSpaced
    \centering
    \caption{Parameters for each case study.}
    \label{fig:results:nj:config}
\begin{tabular}{l|r|r|r}
\hline\hline
Region & New Jersey & Miami & Texas\\
Allocation level & Hospital & Hospital & TSA\\
Bed type & All & All & All\\
Start date & 2020-04-05 & 2020-07-04 & 2020-06-15\\
End date & 2020-06-15 & 2020-08-10 & 2020-08-15\\
Number of days & $72$ & $38$ & $62$\\
Number of locations & $75$ & $24$ & $22$\\
Number of location-days & $5400$ & $912$ & $1364$\\
\end{tabular}
\end{table}
\fi

\iflong
\begin{figure}
    \centering
    \includegraphics[width=0.45\textwidth]{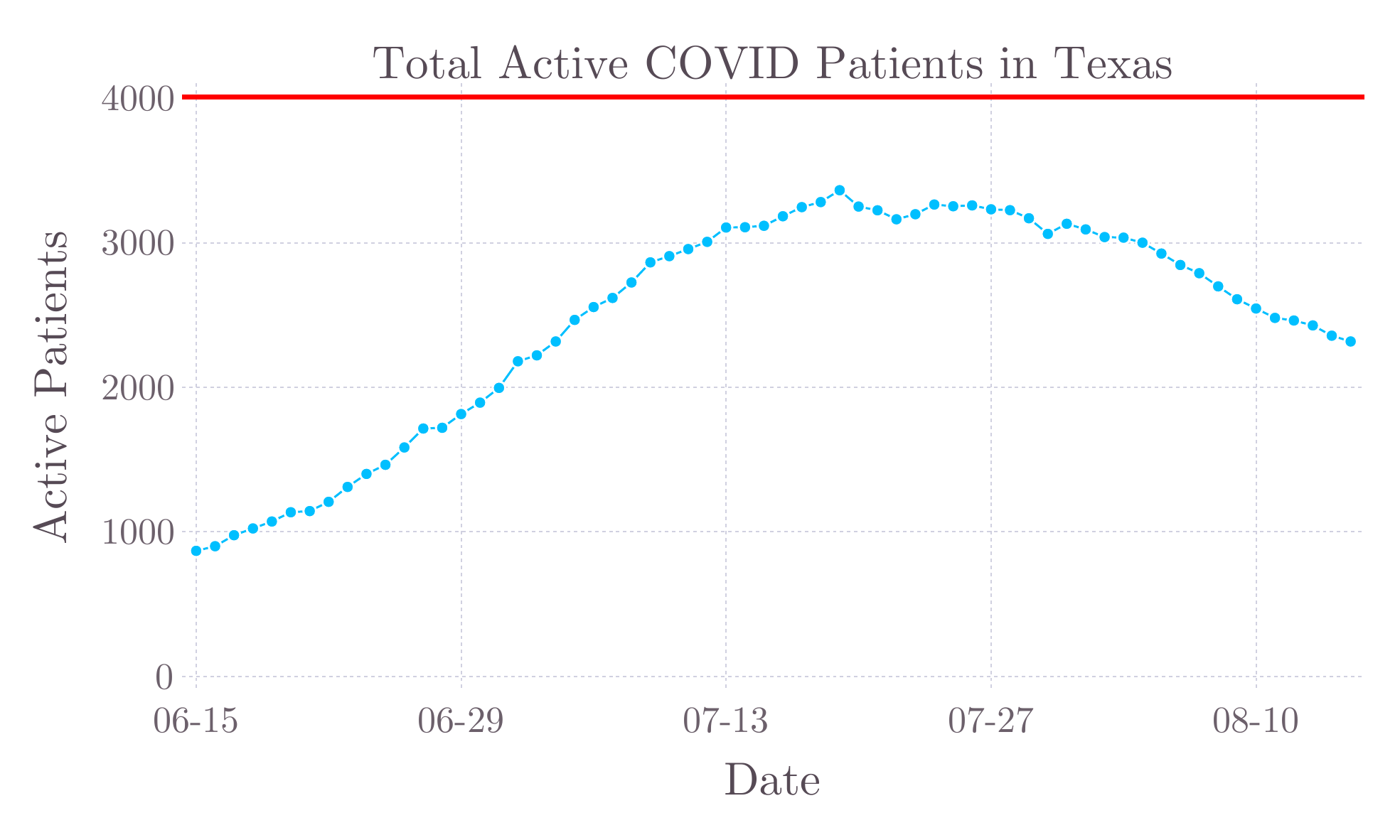}
    \includegraphics[width=0.45\textwidth]{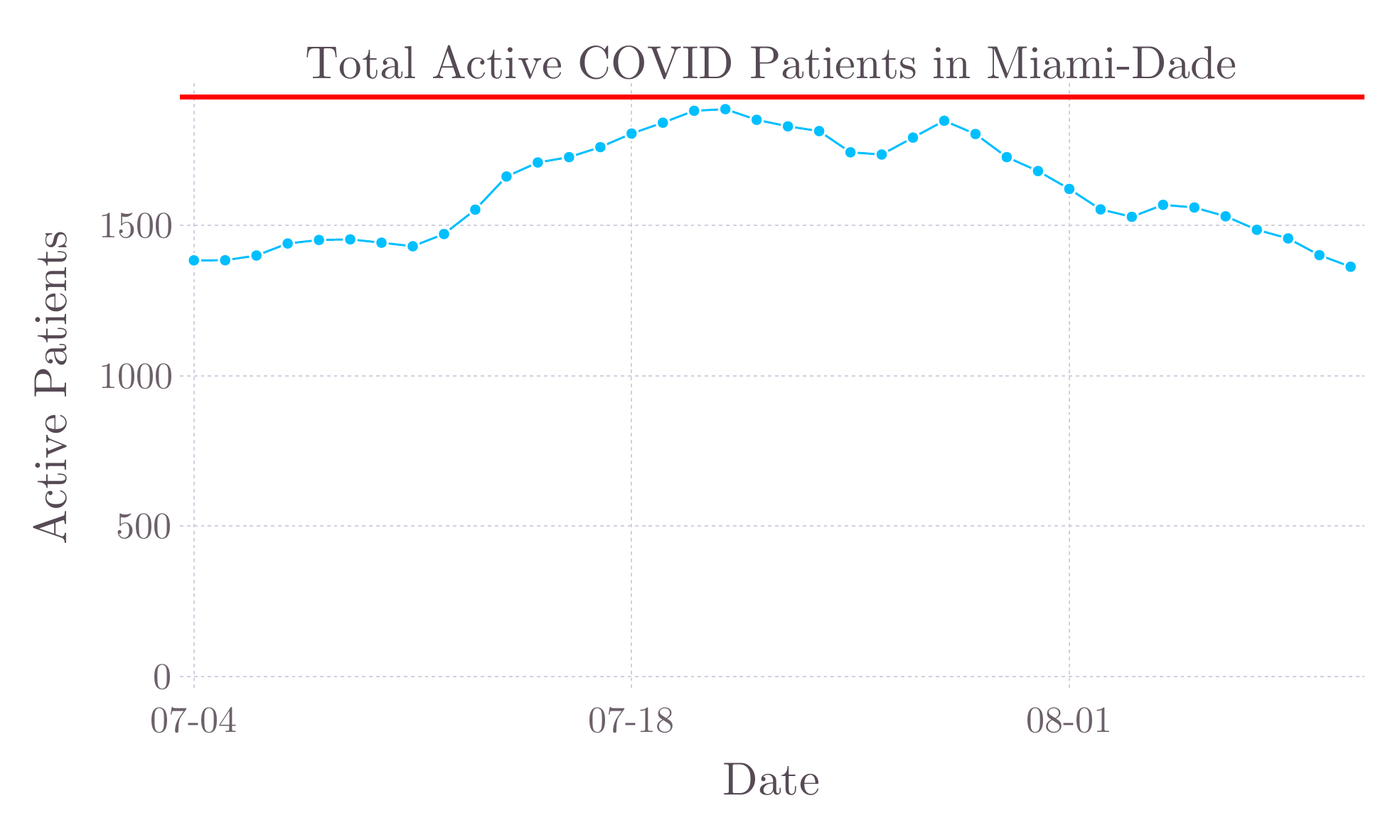}
    \caption{Total active COVID patients (blue) versus the bed capacity for COVID patients (red) for Texas and Miami.}
    \label{fig:results:tx:active_total}
    \label{fig:results:fl:active_total}
    \vspace{-1.5em}
\end{figure}
\fi

\subsubsection{Case Study 1: New Jersey.}
In the previous section we used the New Jersey case study to compare models; here we investigate the solution found by the operational model and its significance for New Jersey. Once again, we consider each of the 75 hospitals in the state of New Jersey from April $5^\text{th}$ to June $15^\text{th}$, 2020.

As noted before, a number of the New Jersey hospitals experienced a COVID patient overflow and as a result, the hospitals took action to increase their original patient capacity to accommodate the additional patients. We estimate that they had to add at least 33,406 bed-days of capacity, or at least a total of 2,093 beds. According to a report from the New Jersey Hospital Association \citep{NJHA2020}, hospitals actually added at least 2,800 ICU beds. Even with the extra capacity, the report also states that the overall load on the system increased from its typical value of 62\% to 82\%, meaning that there was also an increase in demand for nurses. Evidently, the New Jersey hospital system faced an extreme pressure during the first wave of the pandemic.

Despite this extreme pressure, our model found a solution that would have eliminated more than 88\% of the total overflow, as can be seen in Table \ref{fig:results:nj:metrics}. This translates to 3930 surge bed-days and just 1141 additional beds, which is a huge improvement over not transferring patients.
It is clear from these metrics that New Jersey hospitals would not have been forced to create nearly as much surge capacity, putting a far smaller strain on the system, had they made optimal patient transfers.
This solution would have involved transferring 13\% of all COVID patients in New Jersey. While this represents a very large number of people, we believe it is an acceptable trade-off for the potential improvements in patient care, system efficiency, and lower costs to create surge capacity.
Clearly there is a trade-off between the number of transfers and the amount of overflow reduction, which can be made by increasing the penalty on the number of transferred patients in the model. We selected to use a small penalty just to encourage the model not to make unnecessary transfers, but this is potentially an important tool for decision-makers to find operationally feasible solutions.

\begin{figure}[ht]
    \centering
    \includegraphics[width=1.0\textwidth]{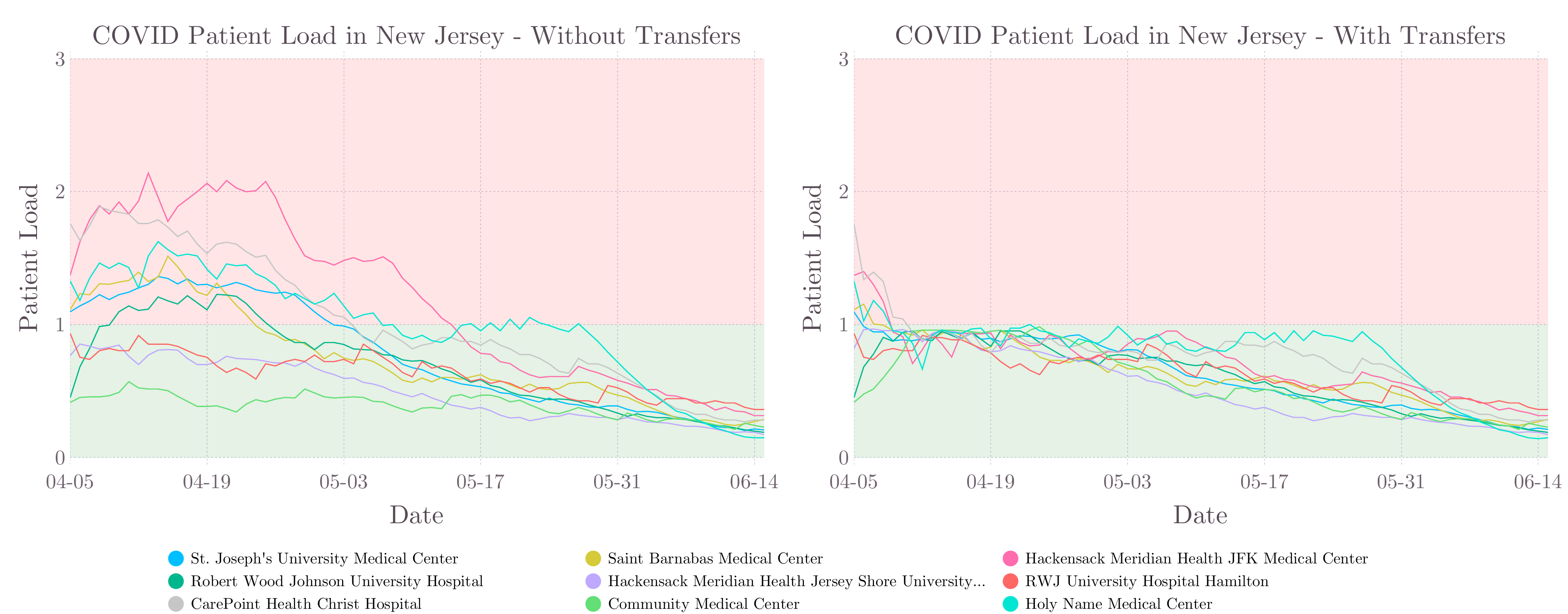}
    \caption{COVID patient load at a representative sample of 9 hospitals in New Jersey, with and without transfers.}
    \label{fig:results:nj:load}
    \vspace{-1.5em}
\end{figure}

\begin{figure}[ht]
    \centering
    \includegraphics[width=0.9\textwidth]{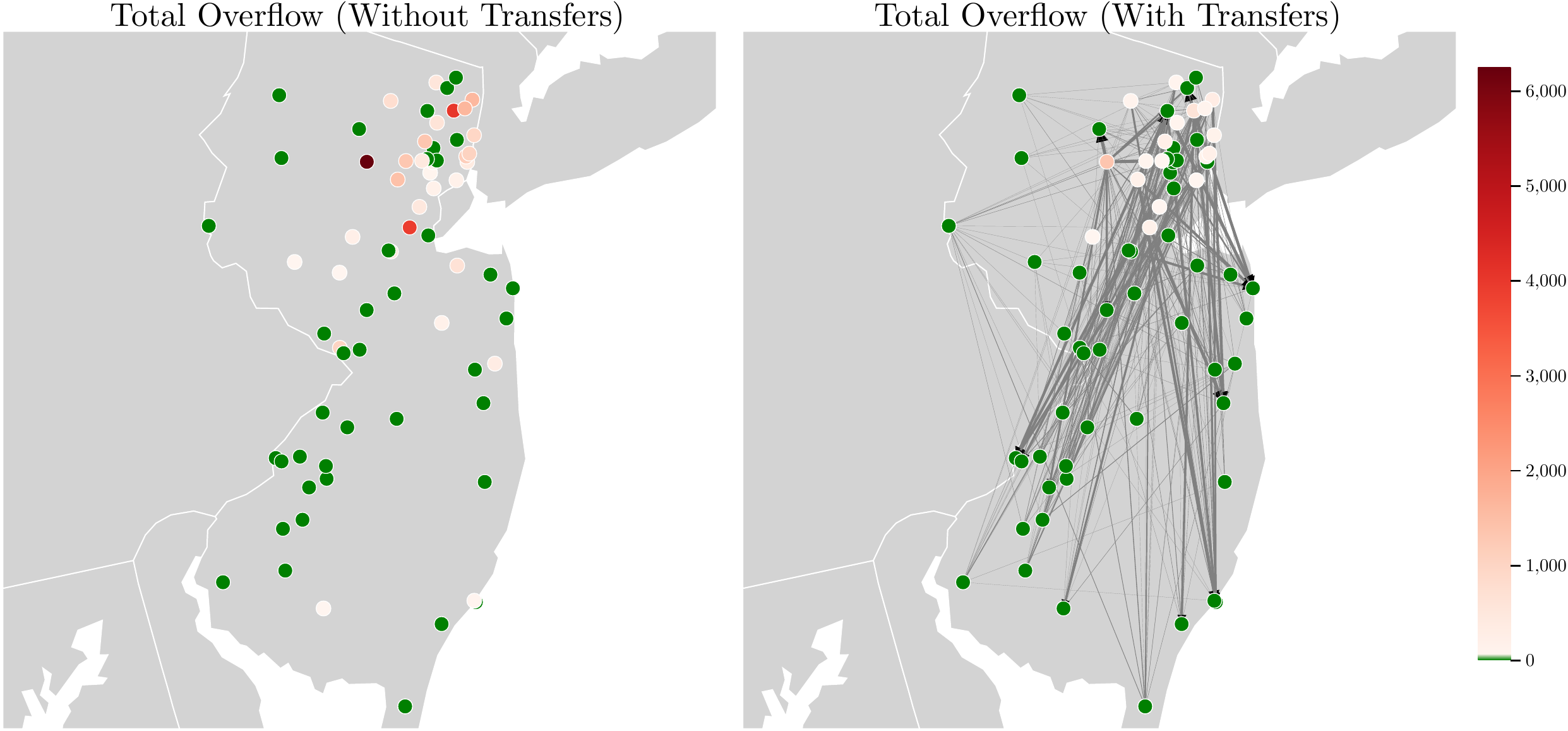}
    \caption{COVID patient overflow by hospital in New Jersey, with and without transfers. Edges indicate patient transfers.}
    \label{fig:results:nj:map}
    \vspace{-1.5em}
\end{figure}

Figure \ref{fig:results:nj:load} plots the COVID patient load over time for a representative sample of hospitals, demonstrating that the model is able to get the additional demand under control almost immediately, and is able to maintain a 5\% buffer between demand and capacity for most hospital-days. In Figure \ref{fig:results:nj:map} we see that the overflow is primarily clustered in the northeast, close to New York City which was the epicenter of the pandemic during this time period, and that the hospitals in the rest of the state were able to stay under capacity. It is this imbalance, and the resulting transfers out of the northeast region, that made the model so successful.

It is also worth noting that while the number of active COVID patients is reported by the hospitals, the exact number of beds the exact surge capacity that could be made available for COVID patients are not known, which adversely affect the accuracy of the results. Accurate collection of such data is paramount to the sound implementation of this system in the real world.

\iflong
\begin{table}
    \TableSpaced
    \centering
    \caption{Evaluation of the performance of the operational model in the New Jersey case study.}
    \label{fig:results:nj:metrics}
    \begin{tabular}{l|r}
\hline\hline
Overflow & $3930$\\
Percent Overflow & $1.2\%$\\
Overflow Reduction & $88.24\%$\\
Percent Of Hospital-Days With An Overflow & $2.41\%$\\
Maximum Overflow By Hospital-Day & $220$\\
Maximum Load By Hospital-Day & $4$\\
\hline
Total Patients Transferred & $3932$\\
Percent Of Patients Transferred & $13.04\%$\\
Max Sent By Day & $483$\\
Max Sent By Hospital & $455$\\
Percent Of Hospital-Days With A Transfer & $24.31\%$\\
Median Non-Zero Transfer & $1$\\
Mean Non-Zero Transfer & $2.4$\\
Max Non-Zero Transfer & $19$\\
\end{tabular}
\end{table}
\fi

\subsubsection{Case Study 2: Texas.}
The second case study considered in this section is the state of Texas at the Trauma Service Area (TSA) level from June to August 2020 when regions of Texas experienced a severe wave of increased COVID cases and some hospitals even practiced transferring patients to balance their load. TSAs are the smallest level at which the state of Texas reports hospitalizations. Analyzing Texas at the TSA level demonstrates that our methods can be applied to healthcare systems at a coarser level while remaining effective and valuable. This approach may be valuable in cases where data is unavailable at more local levels or when entire regions of hospitals are all out or nearly out of capacity. However, such an approach implicitly assumes that there is in fact an optimal distribution of patients among hospitals in each TSA. While this assumption does not hold in general, decision-makers can use the model for each region with particularly unbalanced COVID patient load given the results of the system-wide model.

\begin{figure}[ht]
    \centering
    \includegraphics[width=1.0\textwidth]{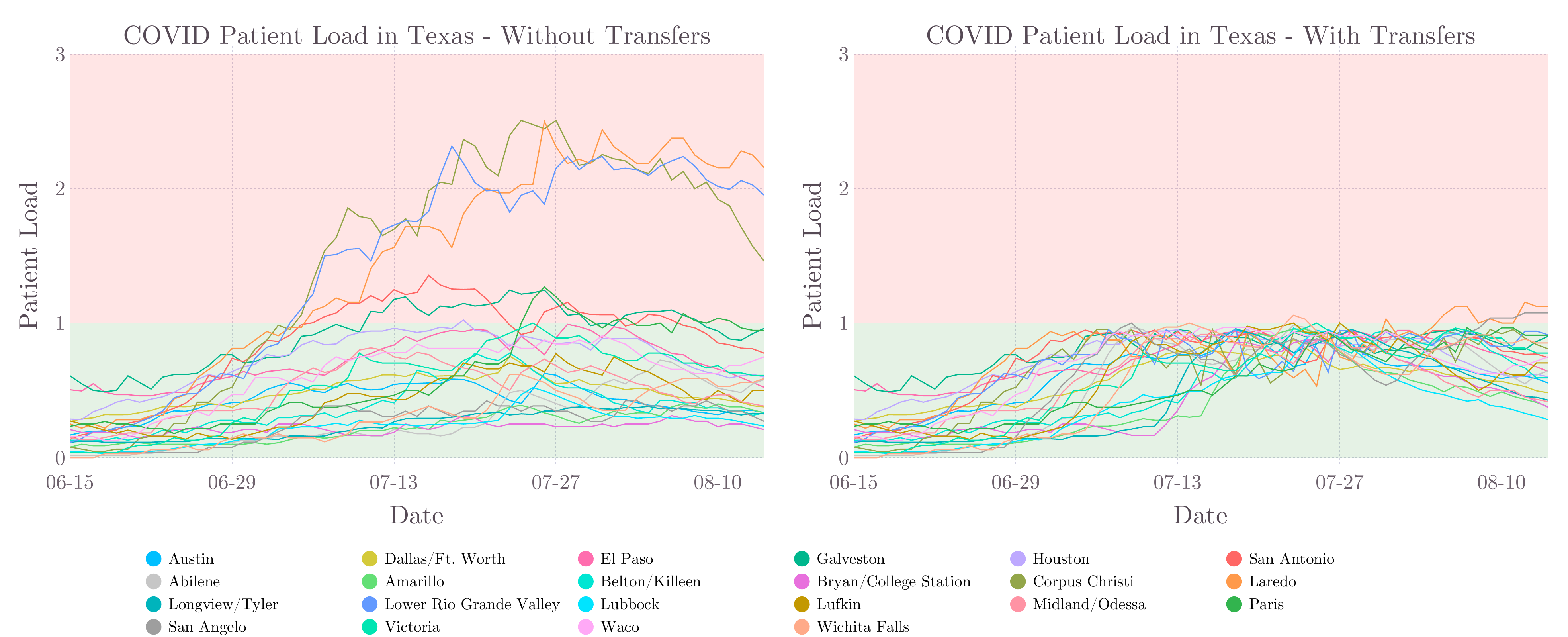}
    \caption{COVID patient load by TSA, with and without transfers.}
    \label{fig:results:tx:load}
    \vspace{-1.5em}
\end{figure}

Figure \ref{fig:results:tx:load} shows that a number of the TSAs go well over a patient load ratio of 1.0, meaning that they had to create a significant amount of surge capacity to care for the COVID patients. However, it should be noted that, similar to the case of the state of new jersey, the surges in the demand were regional and the state did not surpass its capacity as a whole. This, in turn, enabled the models to reduce the total patient overflow to zero by transferring patients among nodes. The results of the transfers can be seen in Table \ref{fig:results:tx:metrics} and Figure \ref{fig:results:tx:load}. The operational model was also able to maintain a gap of 5\% of capacity between the number of active patients and the capacity. The solution of this particular model involved transferring 1765 patients, or 12.1\% of the total COVID patient population for this time period.

\iflong
\begin{table}
    \TableSpaced
    \centering
    \caption{Evaluation of the performance of the operational model in the Texas case study.}
    \label{fig:results:tx:metrics}
    \begin{tabular}{l|r}
\hline\hline
Overflow & $42$\\
Overflow Reduction & $99.7\%$\\
Median Non-Zero Overflow & $2$\\
Mean Non-Zero Overflow & $2.3$\\
Maximum Overflow & $5$\\
Median Load & $70\%$\\
Mean Load & $60\%$\\
Max Load & $120\%$\\
Percent Of Hospital-Days With An Overflow & $1.3\%$\\
\hline
Total Patients Transferred & $1719$\\
Percent Of Patients Transferred & $12.2\%$\\
Percent Of Hospital-Days With A Transfer & $40.9\%$\\
Median Non-Zero Transfer & $2$\\
Mean Non-Zero Transfer & $3.7$\\
Max Non-Zero Transfer & $25$\\
\hline
\end{tabular}
\end{table}
\fi

\begin{figure}[ht]
    \centering
    \includegraphics[width=1.0\textwidth]{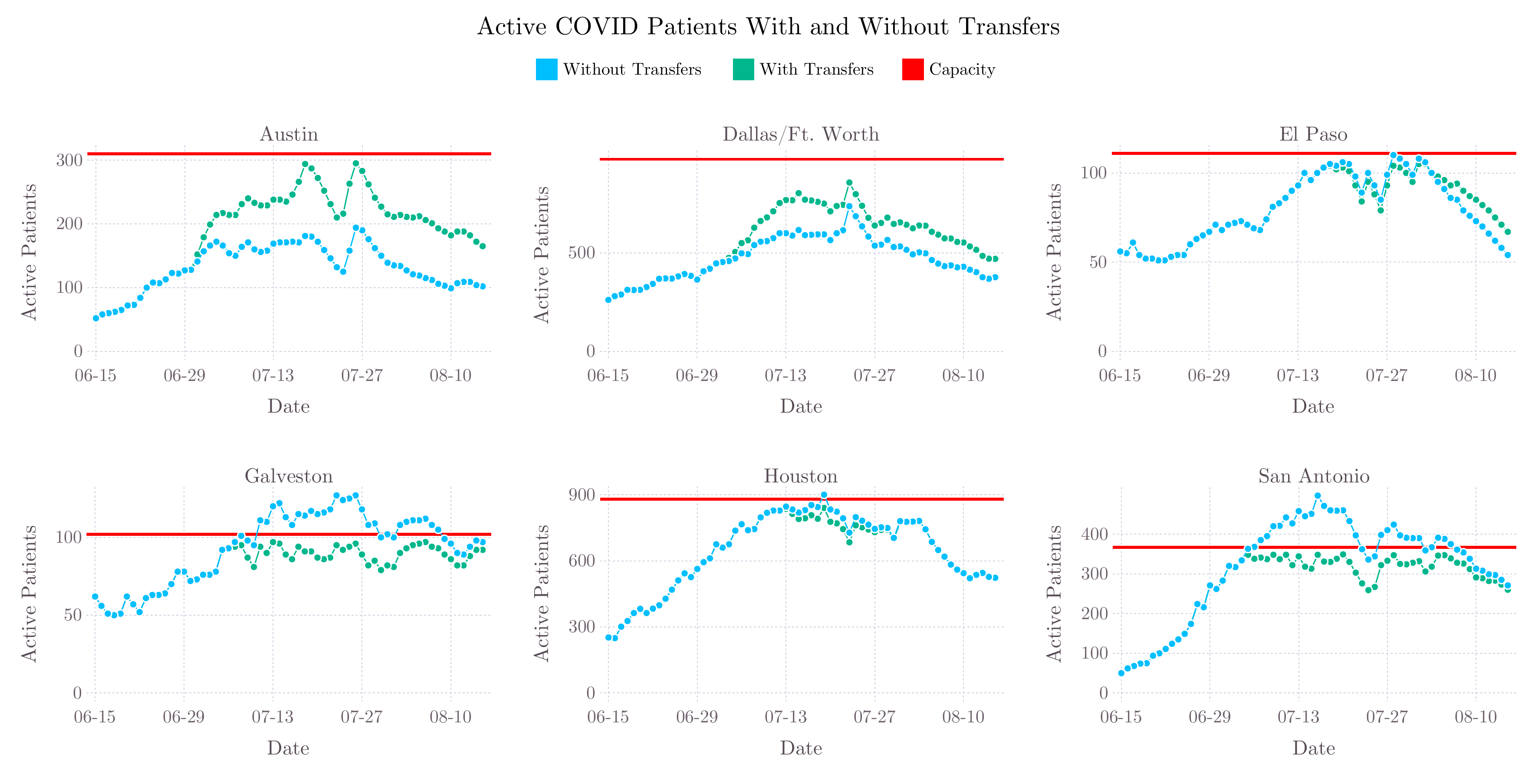}
    \caption{Active COVID patients for a subset of TSAs in Texas, with and without patient transfers.}
    \label{fig:results:tx:active}
    \vspace{-1.5em}
\end{figure}

\iflong
\subsubsection{Case Study 3: Miami}

Our third case study considers hospitals in Florida in July and early August 2020, encompassing the peak of the first wave in the state. We specifically target the 24 Class I hospitals in Miami-Dade County as it was among the most severely impacted counties in Florida at the time. We therefore are investigating the potential value of performing patient re-distribution in a local hospital system rather than state-wide.
We ran the operation model on this system with the same parameters as in the other case studies. The results of this model can be seen in Table \ref{fig:results:fl:metrics}.

According to the evaluation metrics reported in Table \ref{fig:results:fl:metrics}, the operational model could have reduced COVID patient overflow in Miami-Dade by nearly 90\% while requiring fewer than 4\% of the total COVID patients to be transferred. Such a large decrease in overflow and small number of required transfers represents an easy and effective way to accommodate the surge of COVID patients that the system had to face.
The data shows that these hospitals had to increase their capacity by 1251 patients at the peak of the pandemic to accommodate COVID patients. Being able to reduce this by as much as 90\% would have significantly reduced the burden on many hospitals and made ensuring that all patients were properly cared for much easier, cheaper, and more efficient.

\begin{table}
    \TableSpaced
    \centering
    \caption{Evaluation of the performance of the operational model in the Florida case study.}
    \label{fig:results:fl:metrics}
    \begin{tabular}{l|r}
\hline\hline
Overflow & $51$\\
Overflow Reduction & $89.7\%$\\
Median Non-Zero Overflow & $1$\\
Mean Non-Zero Overflow & $1.7$\\
Max Non-Zero Overflow & $21$\\
Median Load & $88.0\%$\\
Mean Load & $85.8\%$\\
Max Load & $121.0\%$\\
Percent Of Hospital-Days With An Overflow & $3.3\%$\\
\hline
Total Patients Transferred & $251$\\
Percent Of Patients Transferred & $3.9\%$\\
Median Non-Zero Transfer & $1$\\
Mean Non-Zero Transfer & $0.6$\\
Max Non-Zero Transfer & $6$\\
Percent Of Hospital-Days With A Transfer & $42.1\%$\\
\end{tabular}
\end{table}

\iflong
\else
\begin{figure}[H]
    \centering
    \includegraphics[width=0.5\textwidth]{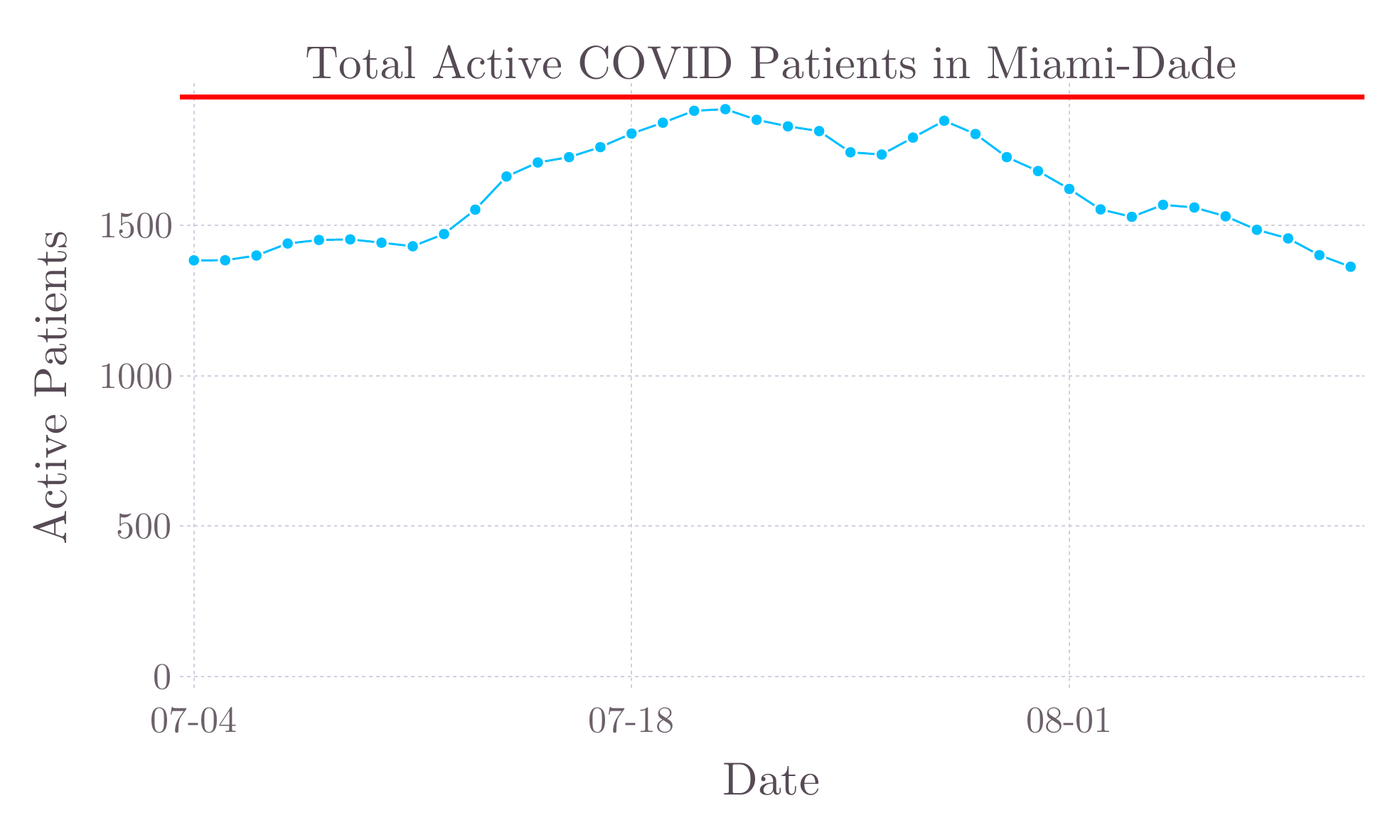}
    \caption{Total active COVID patients (blue) versus the bed capacity for COVID patients (red) in Miami-Dade County.}
    \label{fig:results:fl:active_total}
\end{figure}
\fi

\begin{figure}[H]
    \centering
    \includegraphics[width=1.0\textwidth]{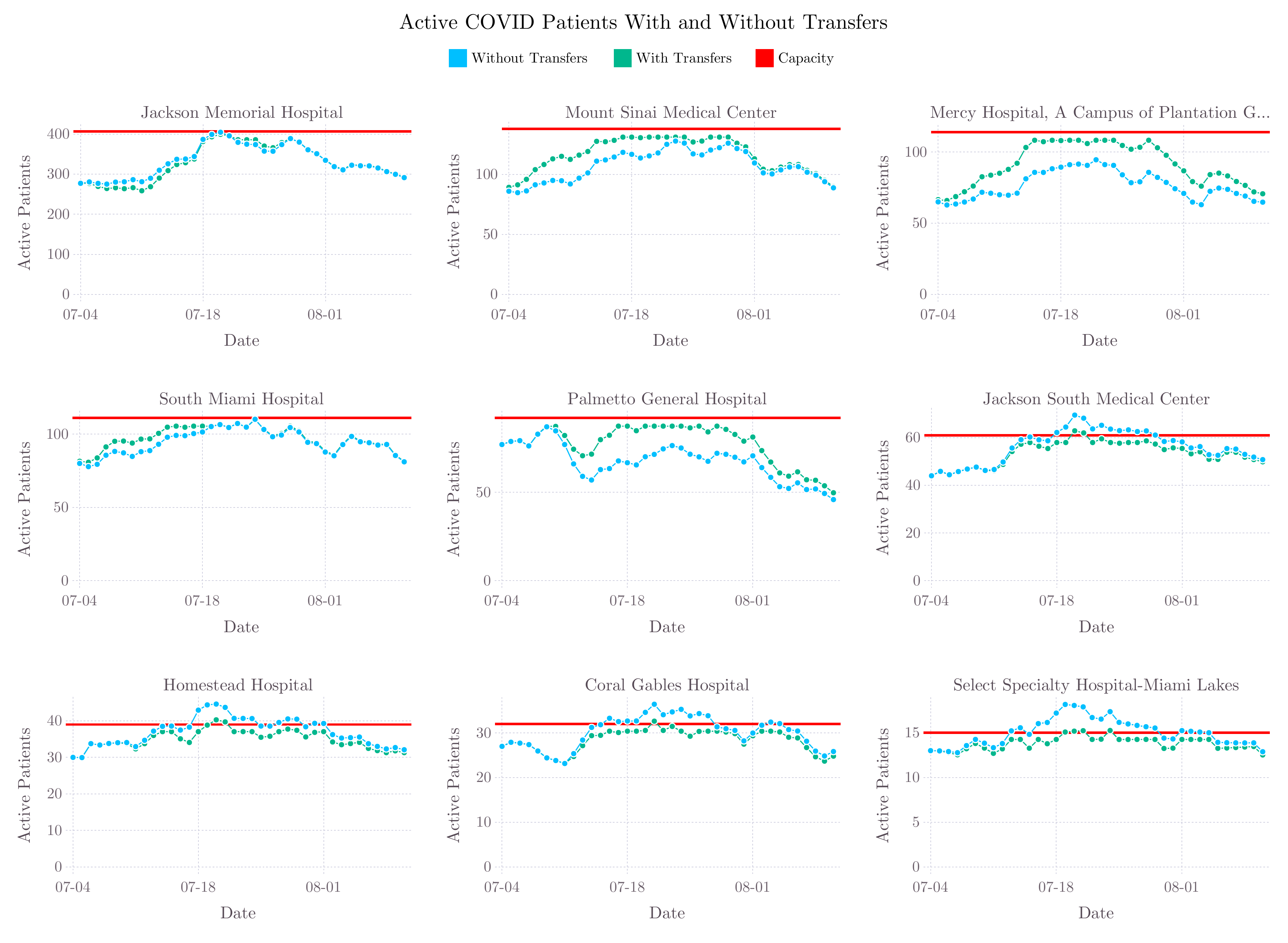}
    \caption{Active COVID patients for 9 representative hospitals in Miami-Dade County, with and without patient transfers.}
    \label{fig:results:fl:active}
\end{figure}

\fi

\subsection{Implementation}
All models described in Section \ref{section:methods} were implemented in Julia 1.5 using the JuMP library for modeling \citep{dunning2017jump}. Models were solved using Gurobi 9.0.3 with default options.
\ifblind
\else
The code is publicly available at: \url{https://github.com/flixpar/covid-resource-allocation}.
\fi

\subsection{Website}
We have demonstrated the potential of our methods to reduce the burden on hospitals facing an extreme surge in demand. In order to make this a practical tool that hospital systems or governments can use we have publicly released the code and data that we have used. In addition to this, we have deployed an interactive website that people can use to explore the potential impact that optimal patient transfers could have. It allows users to specify the region, time period, and model parameters, and run our model in real time. It then displays metrics and animated figures which enables users to effectively understand and evaluate the solution without detailed knowledge of the workings of the model.
\ifblind
\else
The URL for this website is \url{https://covid-hospital-operations.com/}.
\fi

\section{Conclusion} \label{section:conclusion}

In this work, we introduced a methodology that can be used to redistribute demand and resources among hospitals during high demand periods such as the COVID pandemic. The results of applying these methods to real examples from the first wave of the pandemic showed that redistribution effectively helps hospital systems balance patient loads and potentially provide better care while reducing the amount of surge capacity required. Additionally, such a framework promotes systematic collaboration across healthcare entities to proactively respond to predicted high demand events.

The problem of optimizing hospital system decisions during a stressful time like a pandemic is complex and needs proper attention to operational constraints specific to this particular problem. In this work, we introduced a flexible model that can accommodate the varying requirements of hospital systems and optimally allocate newly admitted patients among them to balance patient loads and, consequently, promote patient care quality while reducing operating costs. We introduced different objectives for minimizing the total surge capacity and balancing the load among hospitals, aiming to minimize the total overflow and distributing the stress between the different hospitals. The results of applying these models retrospectively to real-world systems at different levels of the analysis demonstrate successful load balancing, overflow reduction, and good operational characteristics. To improve the practical utility of this work, we also considered resource constraints (specifically the number of available nurses), included different care paths into the model, and made the model robust against some sources of uncertainty. Even with incorporating all such additions and constraints accounting for real-world limiting scenarios and uncertainties, the models outperform no redistribution of patients by a large margin.

Although all the models showcase their capabilities with their exceptional performance in reducing the burden on a healthcare system's entities, there are some limitations to the models that must be addressed. First, primarily while LP formulations (and some MILP) have strong capabilities in efficiently solving large problems enabling us to capture and optimize for healthcare systems on a large scale, but it has its own limitations; it is not possible to accomplish quadratic load balancing or individual patient tracking using such schemes. Additionally, the other limitation is inherent in the data; the scarcity and uncertainty in COVID-related data also pose challenges to using such models. Specifically, due to the decentralized response of the United States to COVID, the quality of infection, hospitalization, personnel, and resource data not only varies widely among counties, states, and regions but also is often poor. Data quality is a serious issue because it is difficult to take decisive action without full information about the situation. However, our robust model targets these issues and alleviates their effects.
It should also be noted that the operational constraints provided do not capture the full range of operational considerations that decision-makers will need to consider. However, the models are flexible to changes and can easily be refined and adapted for each application.

In conclusion, this work results show how patient and resource re-distribution strategies can be immensely effective in reducing the burden on hospitals in extreme situations such as the COVID pandemic. Additionally, our models prove to be capable tools in planning for such strategies, capturing healthcare systems on a large scale, considering multiple operational constraints and sources of uncertainty, and solving in a time-efficient manner while keeping flexible to be easily adapted for each application.

\iflong
\newpage
\fi
\bibliographystyle{informs2014}
\bibliography{covid}

\end{document}